\documentclass[aap,citesort,MSNbibl,dvips]{arximspdf}
\usepackage{graphics}
\doi{10.1214/09-AAP624}
\volume{20}
\issue{2}
\pubyear{2010}
\firstpage{462}
\lastpage{494}

\makeatletter

\DeclareMathAlphabet\mathcaligr{OMS}{cmsy}{m}{n}

\newproclaim{definition}{Definition}[section]
\newtheorem{proposition}{Proposition}[section]
\newtheorem{theorem}[proposition]{Theorem}
\newtheorem{lemma}[proposition]{Lemma}
\newproclaim{remark}{Remark}[section]
\newproclaim{mainquestion}{Main Open Question}

\makeatother

\begin{document}
\begin{frontmatter}

\title{On rough isometries of Poisson processes on~the~line}
\runtitle{Rough isometry of Poisson processes}

\begin{aug}
\author[A]{\fnms{Ron} \snm{Peled}\corref{}\ead[label=e1]{peled@cims.nyu.edu}}
\runauthor{R. Peled}
\affiliation{New York University}
\address[A]{Warren Weaver Hall\\
Courant Institute\\
\quad of Mathematical Sciences\\
New York University\\
251 Mercer Street\\
New York, New York 10012\\
USA\\
\printead{e1}} 
\end{aug}
%

\received{\smonth{9} \syear{2007}}
\revised{\smonth{1} \syear{2010}}

%
\begin{abstract}
Intuitively, two metric spaces are rough isometric (or quasi-isometric)
if their large-scale metric structure is the same, ignoring fine
details. This concept has proven fundamental in the geometric study of
groups. Ab\'ert, and later Szegedy and Benjamini, have posed several
probabilistic questions concerning this concept. In this article, we
consider one of the simplest of these: are two independent Poisson
point processes on the line rough isometric almost surely? Szegedy
conjectured that the answer is positive.

Benjamini proposed to consider a quantitative version which roughly
states the following: given two independent percolations on $\mathbb
{N}$, for
which constants are the first $n$ points of the first percolation rough
isometric to an initial segment of the second, with the first point
mapping to the first point and with probability uniformly bounded from
below? We prove that the original question is equivalent to proving
that absolute constants are possible in this quantitative version. We
then make some progress toward the conjecture by showing that constants
of order $\sqrt{\log n}$ suffice in the quantitative version. This is
the first result to improve upon the trivial construction which has
constants of order $\log n$. Furthermore, the rough isometry we
construct is (weakly) monotone and we include a discussion of monotone
rough isometries, their properties and an interesting lattice structure
inherent in them.
\end{abstract}

%
\begin{keyword}[class=AMS]
\kwd[Primary ]{60D05}
\kwd{60K35}
\kwd[; secondary ]{82B43}.
\end{keyword}
\begin{keyword}
\kwd{Rough isometry}
\kwd{quasi-isometry}
\kwd{Poisson process}
\kwd{percolation}
\kwd{matching}
\kwd{geometry of random sets}
\kwd{one dimension}.
\end{keyword}

\end{frontmatter}

\section{Introduction}\label{sec1}
The concept of \textit{rough isometry} (sometimes also called \textit
{quasi-isometry} or \textit{coarse quasi-isometry}) of two metric
spaces was introduced by Kanai in \cite{K85} and, in the more
restricted setting of groups, by Gromov in \cite{G81}. Informally, two
metric spaces are rough isometric if their metric structure is the same
up to multiplicative and additive constants. This allows stretching and
contracting of distances, as well as having many points of one space
mapped to one point of the other. For example, $R^d$ and $\mathbb
{Z}^d$ are
rough isometric. This concept has proven fundamental in the geometric
study of groups. On the one hand, the rough isometry concept is
stringent enough to preserve some of the metric properties of the
underlying space. On the other hand, it is loose enough to allow for
large equivalence classes of spaces. For example, rough isometry
preserves (under some conditions) geometric properties of the space
such as volume growth and isoperimetric inequalities \cite{K85}. It
preserves analytic properties such as the parabolic Harnack inequality
\cite{D99} (and also \cite{K05}, Section 2.1) and, in a more
probabilistic context, various estimates on transition probabilities of
random walks (heat kernel estimates) are preserved; again, see \cite
{K05}, Section 2, and the references contained therein. Formally, we
have the following.
\begin{definition}
Two metric spaces $X$ and $Y$ are rough isometric (or quasi-isometric)
if there exists a mapping $T\dvtx X\to Y$ and constants $M,D,R\ge0$
such that:
\begin{longlist}
\item any $x_1, x_2\in X$ satisfy
\[
\frac{1}{M}d_X(x_1,x_2) - D\le d_Y(T(x_1), T(x_2)) \le Md_X(x_1,x_2) + D;
\]
\item for any $y\in Y$, there exists some $x\in X$ such that $d_Y(y,
T(x))\le R$.
\end{longlist}
\end{definition}

The first condition ensures that the metric is not distorted too much
multiplicatively or additively; the second condition implies that the
map is close to being onto. On first inspection, it appears that the
definition is not symmetric in $X$ and~$Y$, but one may easily check
that if such a mapping $T\dvtx X\to Y$ exists, then another mapping
$\widetilde{T}\dvtx Y\to X$ also exists, satisfying the same conditions,
with the roles of $X$ and $Y$ interchanged (and with possibly different
constants).

We will sometimes abbreviate ``rough isometric'' to ``r.i.''

In this article, we are concerned with an aspect of the question of how
large the equivalence classes of rough isometric spaces are. We
investigate this question in a probabilistic setting. Specifically Mikl\'os Ab\'ert asked in 2003 \cite{A03} whether, for a finitely generated group,
two infinite clusters of independent edge percolations on its Cayley
graph are rough isometric almost surely (assuming they exist). In this
generality, the question appeared difficult and so Bal\'azs Szegedy suggested
considering whether two site percolations on $\mathbb{Z}^2$ are rough isometric
(disregarding connectivity properties). When this also appeared
difficult, he suggested considering the case of $\mathbb{Z}$. These questions
have since remained open. Independently, and a short time later, the
$\mathbb{Z}
^d$ questions were also raised by Itai Benjamini (following the related
work~\cite{AB03}) who also introduced a quantitative variant. The
one-dimensional question is easily seen to be equivalent to the
following (see Proposition \ref{Poisson_perc_equiv_prop} below): are
two independent Poisson processes on the line (viewed as random metric
spaces with their metric inherited from $\mathbb{R}$) rough isometric a.s.?
Szegedy conjectured a positive answer to this question.

This question is a form of matching problem, but, unlike some other
matching problems in which we wish to minimize some quantity on
average, or to have it bounded for most points, here, we need to
satisfy the rigid constraints of a rough isometry for all points. To
our aid comes the fact that the Poisson processes are infinite and we
may ``start'' constructing the rough isometry at a particularly
convenient location and use the freedom afforded by large constants to
``plan ahead.'' Unfortunately, this article does not settle this
conjecture, but it makes some modest progress. In the next section, we
prove the equivalence of the problem to several other related problems
involving percolations on the integers and on the natural numbers,
including Benjamini's quantitative variant. Our main result is the
construction of a monotone rough isometry with certain properties
giving a first nontrivial upper bound on the quantitative variant.
Section \ref{monotone_rough_isomery_section} presents a discussion of
monotone rough isometries, their properties and an interesting lattice
structure inherent in them. As noted there, in general, monotone rough
isometries between subsets of $\mathbb{Z}$ are more restrictive than general
rough isometries. In particular, it may be harder to find a monotone
rough isometry between two independent Poisson processes than to find a
general rough isometry. Section \ref{equivalence_thms_proofs_section}
contains the proofs of all the theorems in Section \ref
{variants_section}, except for the main construction. Section \ref
{main_construction_section} presents the main construction.

\section{Versions of the problem and main result} \label{variants_section}
In this section, we will first state in precise terms the main open
question described in the \hyperref[sec1]{Introduction}. We will then
proceed to show
the equivalence of the question to several other related problems. We
shall go from the continuous Poisson process question to a discrete
variant (percolation on $\mathbb{Z}$), then to an oriented discrete variant
(percolation on $\mathbb{N}$) and, finally, to a finite variant
(percolation on
an initial segment of $\mathbb{N}$), all of which are equivalent. We
will then
state a quantitative version of our main open question (due to
Benjamini), based on the finite variant, and conclude the section with
a statement of our main result which gives the first nontrivial upper
bound on this quantitative version. The proofs of all statements in
this section, except for the main result, are presented in Section \ref
{equivalence_thms_proofs_section}; the proof of the main result is
presented in Section \ref{main_construction_section}.
\begin{proposition} \label{zero_one_event_prop}
Given two independent Poisson processes $A,B\subseteq\mathbb{R}$
(possibly of
different intensities) and constants $(M,D,R)$, the event that $A$ and
$B$ are rough isometric with constants $(M,D,R)$ is a zero-one event.
\end{proposition}

Hence, we come to the
following question.
\begin{mainquestion}
Do there exist constants $(M,D,R)$ for which two independent Poisson
processes of intensity 1 are rough isometric a.s.\textup{?}
\end{mainquestion}

In this article, we shall mostly consider a discrete variant of the
question involving Bernoulli percolations on $\mathbb{Z}$ or on
$\mathbb{N}$, rather
than Poisson processes. We remind the reader that a Bernoulli
percolation on $\mathbb{Z}$ with parameter $p$ is the random subset
$A\subseteq
\mathbb{Z}$ obtained from $\mathbb{Z}$ by independently deleting each
integer with
probability $1-p$. It is defined analogously for $\mathbb{N}$. The next
proposition states the equivalence of the problem for Bernoulli
percolations and for Poisson processes.
\begin{proposition} \label{Poisson_perc_equiv_prop}
The following are equivalent:
\begin{longlist}
\item for some intensities $\alpha, \beta>0$, two independent Poisson
processes, one with intensity $\alpha$ and the other with intensity
$\beta$, are rough isometric a.s.;
\item for any intensities $\alpha, \beta>0$, two independent Poisson
processes, one with intensity $\alpha$ and the other with intensity
$\beta$, are rough isometric a.s.;
\item for some $0<p,q<1$, two independent Bernoulli percolations on
$\mathbb{Z}
$, one with parameter $p$ and the other with parameter $q$, are rough
isometric a.s.;
\item for any $0<p,q<1$, two independent Bernoulli percolations on
$\mathbb{Z}
$, one with parameter $p$ and the other with parameter $q$, are rough
isometric a.s.
\end{longlist}
\end{proposition}

Since, by the previous proposition, we may equivalently consider any
intensity for the Poisson process and any parameter for Bernoulli
percolation, we fix notation and, from this point on, consider only
Poisson processes with unit intensity and Bernoulli percolations with
parameter $\frac{1}{2}$.

A rough isometry between two Poisson processes or between two Bernoulli
percolations on $\mathbb{Z}$ is not necessarily order preserving (or order
reversing), as will be discussed in more detail near the end of this
section. Still, one feels intuitively that such a mapping should be
monotonic in some rough sense. Indeed, for the next two theorems, we
will need to show that such a mapping is at least ``almost monotonic at
most points,'' in a sense made precise in the following statements and
their proofs. We start by showing that a certain oriented version of
the problem is equivalent to the original problem. For this purpose, we
introduce the following new concept.
%
\begin{definition}
Two rooted metric spaces $(X,a)$ and $(Y,b)$ are \textit{rooted rough
isometric} if there exists a mapping $T\dvtx X\to Y$ and constants
$M,D,R\ge0$ such that $T(a)=b$ and the conditions in the usual
definition of rough isometry hold for $T$ and the constants $(M,D,R)$.
\end{definition}

We also introduce a different random model, as follows.
\begin{definition}
A \textit{rooted Bernoulli percolation on $\mathbb{N}$} (with parameter
$\frac
{1}{2}$) is a random subset $A\subseteq\mathbb{N}\cup\{0\}$ in which
$0\in A$
deterministically and any $n\in\mathbb{N}$ belongs to $A$ with probability
$\frac{1}{2}$ independently.
\end{definition}
\begin{theorem} \label{discrete_to_oriented_equiv_thm}
The following are equivalent:
\begin{longlist}
\item two independent Bernoulli percolations on $\mathbb{Z}$ (with parameter
$\frac{1}{2}$) are rough isometric a.s.;
\item two independent rooted Bernoulli percolations $(A,0)$ and $(B,0)$
on $\mathbb{N}$ are rooted rough isometric with positive probability.
\end{longlist}
\end{theorem}

To prove this theorem, we need the following
definition.
\begin{definition} \label{cut_point_definition}
Given two subsets $A,B\subseteq\mathbb{Z}$ and a mapping $T\dvtx A\to
B$, the
point $x\in A$ is called a \textit{cut point} for $T$ if one of the
following occurs:
\begin{enumerate}[($\alpha$)]
\item[($\alpha$)] for all $z\in A$ with $z>x$, we have $T(z)\ge T(x)$;
\item[($\beta$)] for all $z\in A$ with $z>x$, we have $T(z)\le T(x)$.
\end{enumerate}
\end{definition}

We also require the following lemma.
\begin{lemma} \label{cut_point_lemma}
If two independent Bernoulli percolations on $\mathbb{Z}$ are rough isometric
a.s. with constants $(M,D,R)$, then, with probability 1, any rough
isometry $T\dvtx A\to B$ with constants $(M,D,R)$ has a cut point.
\end{lemma}

We continue to construct a finite variant of our problem. First, we
define, for a given infinite subset $A\subseteq\mathbb{N}\cup\{0\}$,
$A(n)\subseteq A$ to be its first $n$ points [e.g., if
$A=(0,1,3,4,6,\ldots)$, then $A(3)=(0,1,3)$]. Also, given
$A,B\subseteq
\mathbb{N}\cup\{0\}$ both containing $0$, we sometimes say that
$A(n)$ is
\textit{rooted r.i. to some initial segment of $B$} if there exists an
$m$ and a rooted r.i. $T\dvtx A(n)\to B(m)$. We may also phrase this as
\textit{$T$ is a rooted r.i. of $A(n)$ to some initial segment of $B$}.
We now have the following result.
\begin{theorem} \label{finite_and_infinite_perc_equiv_thm}
The following are equivalent:
\begin{longlist}
\item two independent rooted Bernoulli percolations $(A,0)$ and $(B,0)$
on $\mathbb{N}$ are rooted rough isometric with positive probability;
\item there exists some $p>0$ and constants $(M,D,R)$ such that, given
two independent rooted Bernoulli percolations $(A,0)$ and $(B,0)$ on
$\mathbb{N}
$, for any $n\ge1$, $A(n)$ is rooted r.i. to some initial segment of
$B$ with constants $(M,D,R)$ and with probability at least $p$.
\end{longlist}
\end{theorem}

Although this theorem may initially seem straightforward, it transpires
that the direction $\mbox{(i)} \to \mbox{(ii)}$ is somewhat problematic. The main
difficulty stems from the fact that, given a rooted rough isometry from
$A\subseteq\mathbb{N}\cup\{0\}$ to $B\subseteq\mathbb{N}\cup\{0\}
$, its
restriction to
$A(n)$ is not necessarily a rooted rough isometry to some $B(m)$ with
the same constants. This is due to the fact that a rough isometry need
not be monotonic and hence the image of its restriction to $A(n)$ may
still have big ``holes'' [i.e., points $b\in B$ where property (ii) in
the definition of rough isometry does not hold] which are ``filled'' by
the mapping at subsequent points of $A$. To prove this theorem, we will
need a statement asserting that if $A$ and $B$ are rooted Bernoulli
percolations, then $T\dvtx A\to B$ is a rooted rough isometry with
constants $(M,D,R)$, and if we allow the constant $R$ to be increased
sufficiently, say to some $L:=L(M,D,R)$, then for ``most $n$'s'' the
restriction of $T$ to $A(n)$ will still be a rooted rough isometry to
some initial segment of $B$ with constants $(M,D,L)$. This is the
content of the next three lemmas: they make precise what was meant when
we stated previously that a rough isometry is ``almost monotonic at most
points.''

We now introduce the following notation: given $A\subseteq\mathbb
{N}\cup\{0\}$
and $x\in A$, let $\operatorname{Succ}(x)$ be the smallest point in
$A$ which is
larger than $x$ (or $\infty$ if there is no such point) and let
$\operatorname{Gap}
(x):=\operatorname{Succ}(x)-x$. We start with a deterministic lemma.
\begin{lemma} \label{big_gap_deterministic_lemma}
Let $A,B\subseteq\mathbb{N}\cup\{0\}$ be infinite subsets, both
containing $0$,
and let $T\dvtx A\to B$ be a rooted r.i. between them with constants
$(M,D,R)$. There exists $L:=L(M,D)$ such that if there exist $x,y\in A$
with $x<y$ and $T(y)\le T(x)-L$, then there exists $z\in A$, $z\ge y$
and $z-x\ge\frac{L}{2M}$ such that $\operatorname{Gap}(z)\ge\frac
{z-x}{2M^2}$.
\end{lemma}

We continue with a probabilistic aspect of the previous lemma.
\begin{lemma} \label{big_gap_prob_estimate_lemma}
Let $A$ be a rooted Bernoulli percolation on $\mathbb{N}$, let $w\in
\mathbb{N}$ and
define, for constants $L,M$, the event $E^w_{L,M}:=\{\exists z\in A,
z>w, \operatorname{Gap}(z)\ge\max(\frac{L}{4M^3},\frac{z-w}{2M^2}
)\}$. Then $\mathbb{P}
(E^w_{L,M})\le C(\frac{L}{M}+1)e^{-c{L}/{M^3}}$ for some absolute
constants $C,c>0$ (not depending on any parameter).
\end{lemma}

Finally, we have one more deterministic lemma.
\begin{lemma} \label{rooted_ri_restriction_lemma}
Let $A,B\subseteq\mathbb{N}\cup\{0\}$ be infinite subsets, both
containing $0$,
and let $T\dvtx A\to B$ be a rooted r.i. between them with constants
$(M,D,R)$. Fix $L>R$ and $n\ge1$, let $x_n$ be the $n$th point of $A$
and suppose that the event $E^{x_n}_{L-R,M}$ of Lemma~\ref
{big_gap_prob_estimate_lemma} does not hold for $A$. Then $T$
restricted to $A(n)$ is a rooted r.i. of $A(n)$ to $B(m)$ for some $m$
with constants $(M,D,L)$.
\end{lemma}
\begin{remark}
Close inspection of the proof of part $\mbox{(ii)}\to\mbox{(i)}$ of Theorem \ref
{finite_and_infinite_perc_equiv_thm} reveals that (ii) is, in fact,
equivalent (by the same proof) to the following, seemingly weaker,
statement [the $R$-denseness property is property (ii) in the
definition of r.i.]:

(iii) There exists $p>0$, constants $(M,D,R)$ and a function $f(n)\to
\infty$ such that given two independent rooted Bernoulli percolations
$(A,0)$ and $(B,0)$ on $\mathbb{N}$, for any $n\ge1$, with
probability at
least $p$, there exists $T$ from $A(n)$ to $B(m)$ for some $m$ (a
function of $A$, $B$ and $n$) which satisfies the properties of a
rooted rough isometry with constants $(M,D,R)$, except that we only
require the $R$-denseness property to hold for $b\in B(m)$ with $b\le f(n)$.

Since this statement is complicated to state and we make no use of it
in the sequel, we simply leave it as a remark.
\end{remark}

The last theorem gives rise to the following quantitative variant of
our main question which will be our main concern in this article.
\begin{mainquestion} \label{quantitative_open_question}
Given two independent rooted Bernoulli percolations $(A,0)$ and $(B,0)$
on $\mathbb{N}$, for which functions $(M(n), D(n), R(n))$ does there
exist a
rooted rough isometry $T$ with constants $(M(n),D(n),R(n))$ from
$(A(n),0)$ to $(B(m),0)$ for some $m$ (a function of $A$, $B$ and $n$)
with probability not tending to 0 with $n$\textup{?}
\end{mainquestion}

By the previous theorem, our first open question is equivalent to the
claim that constant functions suffice. Both the original open question
and this quantitative variant were posed to the author by Itai
Benjamini \cite{B05} (although without the proof of equivalence) and it
is the main aim of this paper to present some progress on this
quantitative variant.

Trivially, one has that the functions $(\log_2 n,0,0)$, or even $(\log
_2 n - C, 0, 0)$ for some $C>0$, suffice for this quantitative question
by considering the mapping from $(A(n),0)$ to $(B(n),0)$ which maps the
$i$th point of $A$ to the $i$th point of $B$. We are not aware of any
improvement on this trivial result in the literature. We can now state
our main result.
\begin{theorem} \label{main_thm}
There exists $N>0$ such that, given two independent rooted Bernoulli
percolations $(A,0)$ and $(B,0)$ on $\mathbb{N}$, for any $n>N$, there
exists a
random $m$ (a function of $A$, $B$ and $n$) such that $(A(n),0)$ and
$(B(m),0)$ are rooted rough isometric with constants $(30\sqrt{\log_2
n},\frac{1}{2},10\sqrt{\log_2 n})$ and with probability
$1-2^{-8\sqrt
{\log_2 n}}$.
\end{theorem}

This theorem is proved by a direct construction which will be detailed
in Section \ref{main_construction_section}. Furthermore, the mapping we
construct is (weakly) monotone increasing (in fact, we construct a
Markov rough isometry in the sense of Section \ref
{Markov_ri_subsection}). As already noted, monotonicity is not required
by the definition of rough isometry, but monotone mappings are easier
to construct, have nicer properties and an interesting structure, as
explained in the next section. We do not know if the question of having
a monotone rough isometry between (say) two Poisson processes is
equivalent to the question of having just a general rough isometry
between them. The next section also makes this question precise.
\begin{remark}
We note that, up to the constant $8$, the success probability achieved
in Theorem \ref{main_thm} is optimal. To see this, consider the event $(0,1,\ldots,
\lceil15\sqrt{\log_2 n}\rceil+ 1)\subseteq A$ and that in $B$ the
next point after $0$ is greater than $30\sqrt{\log_2 n} + \frac{1}{2}$.
This event has probability larger than $2^{-45\sqrt{\log_2 n}-3}$ and
we claim that on it there is no rooted r.i. between $A$ and $B$
with constants $(30\sqrt{\log_2 n},\frac{1}{2},10\sqrt{\log_2
n})$. To
see this, suppose, in order to reach a contradiction, that there was
such a rooted r.i.~$T$. If we let $x_0=\max(x\in A | T(x)=0)$, then we
must have $x_0\le15\sqrt{\log_2 n}$ by property (i) of the r.i. [and
since $T(0)=0$]. Hence, $x_0+1\in A$ and we must have $T(x_0+1)>30\sqrt
{\log_2 n} + \frac{1}{2}$. This is a contradiction since, then,
$30\sqrt
{\log_2 n} + \frac{1}{2}<T(x_0+1) - T(x_0)\le30\sqrt{\log_2
n}(x_0+1 -
x_0)+\frac{1}{2}$.
\end{remark}

\section{Monotone rough isometries}\label{monotone_rough_isomery_section}
\begin{figure}

\includegraphics{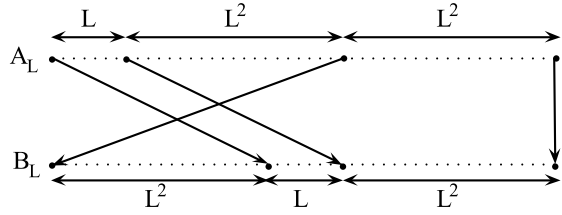}

\caption{Monotonic rough isometry must have large constants.}
\label{monotone_non_monotone_L_example}
\end{figure}

\begin{figure}[b]

\includegraphics{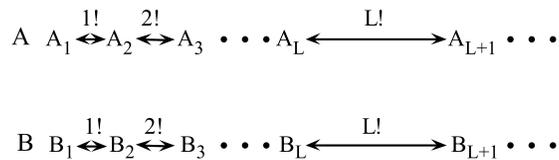}

\caption{No monotonic rough isometry exists.}
\label{monotone_non_monotone_L_many_example}
\end{figure}

In this section, we consider the notion of a (\textit{weakly}) \textit
{increasing rough isometry}, that is, a rough isometry mapping $T\dvtx
X\to Y$ between two subsets $X,Y\subseteq\mathbb{R}$ for which
$T(x)\ge T(y)$
whenever $x\ge y$. As is easy to check, the notion of an increasing
rough isometry defines an equivalence class on subsets of $\mathbb
{R}$, that
is, if $X,Y,Z\subseteq\mathbb{R}$ and $T_1\dvtx X\to Y$, $T_2\dvtx
Y\to Z$ are
increasing rough isometries, then there also exist $T_3\dvtx Y\to X$
and $T_4\dvtx X\to Z$ ($T_4:=T_2\circ T_1$) which are increasing rough
isometries. If there exists an increasing rough isometry between such
$X$ and $Y$, we shall call $X$ and $Y$ \textit{increasing rough
isometric}. On first reflection, one may hope that the notions of
increasing rough isometry and general rough isometry are equivalent,
that is, that if two spaces $X,Y\subseteq\mathbb{R}$ are rough
isometric, then
they are also increasing rough isometric (perhaps with different
constants). Unfortunately, this is not the case in general, as one may
see by means of various examples. Figures \ref
{monotone_non_monotone_L_example} and \ref
{monotone_non_monotone_L_many_example} show a variant of an example
shown to the author by Gady Kozma \cite{K06}. For each integer $L\ge
1$, Figure \ref{monotone_non_monotone_L_example} shows two subsets
$A_L,B_L\subseteq\mathbb{N}$ (each containing four points), between
which there
exists a nonmonotone rough isometry with constants $(3,0,0)$ (which is
depicted). However, as is easy to see, any (weakly) monotone rough
isometry will have constants tending to infinity with the parameter $L$.

%

%

Although this example involves two finite sets of points and, of
course, any two finite sets are increasing rough isometric for some
constants, one may use this example to construct two infinite sets of
points which are rough isometric, but not increasing rough isometric.
Figure \ref{monotone_non_monotone_L_many_example} shows two such sets
$A$ and $B$ which are constructed by concatenating the previous $(A_L,
B_L)$ example, but with a gap of size $L!$ in both $A$ and $B$ between
$(A_L,B_L)$ and $(A_{L+1},B_{L+1})$. On the one hand, concatenating the
rough isometries of Figure \ref{monotone_non_monotone_L_example} gives
a rough isometry with finite constants here, but, on the other hand,
such a fast growing gap ensures that any rough isometry between $A$ and
$B$ will have some large $L$ (depending on its constants) such that for
all $j\ge L$, the points of $A_j$ will only be mapped to the points of
$B_j$, thereby reducing to the example of Figure \ref
{monotone_non_monotone_L_example} within each such segment. In
particular, the rough isometry cannot be monotonic.

In our context, it is then natural to ask
the following question.
\begin{mainquestion} \label{Monotone_Poisson_RI_question}
Given two independent Poisson processes $A,B$, does there exist a
(weakly) increasing rough isometry between them a.s.\textup{?}
\end{mainquestion}

As in Section \ref{variants_section}, one can prove the following.
\begin{proposition}
Given two independent Poisson processes $A,B\subseteq\mathbb{R}$
(possibly of
different intensities) and constants $(M,D,R)$, the event that $A$ and
$B$ are increasing rough isometric with constants $(M,D,R)$ is a
zero-one event.
\end{proposition}

We also have the following equivalences.
\begin{proposition}
The following are equivalent:
\begin{longlist}
\item for some intensities $\alpha, \beta>0$, two independent Poisson
processes, one with intensity $\alpha$ and the other with intensity
$\beta$, are increasing rough isometric a.s.;
\item for any intensities $\alpha, \beta>0$, two independent Poisson
processes, one with intensity $\alpha$ and the other with intensity
$\beta$, are increasing rough isometric a.s.;
\item for some $0<p,q<1$, two independent Bernoulli percolations on
$\mathbb{Z}
$, one with parameter $p$ and the other with parameter $q$, are
increasing rough isometric a.s.;
\item for any $0<p,q<1$, two independent Bernoulli percolations on
$\mathbb{Z}
$, one with parameter $p$ and the other with parameter $q$, are
increasing rough isometric a.s.
\end{longlist}
\end{proposition}

The proofs of these statements are exactly the same as in Section \ref
{variants_section}, but with ``rough isometry'' replaced by
``increasing rough isometry,'' and are hence omitted. Again, due to
these equivalences, we shall only consider Poisson processes of unit
intensity and Bernoulli percolations with parameter $\frac{1}{2}$.

Analogously to Section \ref{variants_section}, we can define a \textit
{rooted increasing rough isometry} between two rooted spaces $(X,a)$
and $(Y,b)$, where $X,Y\subseteq\mathbb{R}$, as a mapping $T\dvtx
X\to Y$ which
is an increasing rough isometry and has $T(a)=b$. For increasing rough
isometries, it is much easier to pass from the question about
percolations on $\mathbb{Z}$ to the question about percolations on
$\mathbb{N}$, and
from there to the finite version. This is due to the following
obvious statement.
\begin{proposition} \label{increasing_RI_restriction_prop}
If $A,B\subseteq\mathbb{R}$ are increasing rough isometric by a mapping
$T\dvtx
A\to B$ with constants $(M,D,R)$, then, for any $x,y\in A$ with $x<y$,
we have that $T$ restricted to $A\cap[x,y]$ is an increasing rough
isometry from $A\cap[x,y]$ to $B\cap[T(x),T(y)]$ with constants $(M,D,R)$.
\end{proposition}

We emphasize once more that this statement is not true for general
rough isometries, although, for increasing rough isometries, it is
trivial to check that it holds (we omit the proof). From this, we
easily deduce the following result.
\begin{theorem}
The following are equivalent:
\begin{longlist}
\item two independent Bernoulli percolations on $\mathbb{Z}$ are increasing
rough isometric a.s.;
\item two independent rooted Bernoulli percolations $(A,0)$ and $(B,0)$
on $\mathbb{N}$ are rooted increasing rough isometric with positive
probability;
\item there exists some $p>0$ and constants $(M,D,R)$ such that, given
two independent rooted Bernoulli percolations $(A,0)$ and $(B,0)$ on
$\mathbb{N}
$, for any $n\ge1$, $A(n)$ is a rooted increasing r.i. to some initial
segment of $B$ with constants $(M,D,R)$ and with probability at least $p$.
\end{longlist}
\end{theorem}

Using Proposition \ref{increasing_RI_restriction_prop}, the
equivalences $\mbox{(i)}\to\mbox{(ii)}$ and $\mbox{(ii)}\to\mbox{(iii)}$ are trivial to prove. The
proofs of $\mbox{(ii)}\to\mbox{(i)}$ and $\mbox{(iii)}\to\mbox{(ii)}$ are the same as those given
in Theorems \ref{discrete_to_oriented_equiv_thm} and \ref
{finite_and_infinite_perc_equiv_thm}, with ``rough isometry'' replaced
by ``increasing rough
isometry.''

Of course, one can now formulate a quantitative version of our
question, as follows.
\begin{mainquestion} \label{increasing_quantitative_open_question}
Given two independent rooted Bernoulli percolations $(A,0)$ and $(B,0)$
on $\mathbb{N}$, for which functions $(M(n), D(n), R(n))$ does there
exist an
increasing rooted rough isometry $T$ with constants $(M(n),D(n),\break R(n))$
from $(A(n),0)$ to $(B(m),0)$ for some $m$ (a function of $A$, $B$ and
$n$) with probability not tending to 0 with $n$\textup{?}
\end{mainquestion}

As was mentioned earlier, Theorem \ref{main_thm} is still relevant in
this context since the rough isometries we construct there are
increasing rough isometries.

Until now, we have stated the common features of general rough
isometries and increasing rough isometries. The next two subsections
present some features which are unique to increasing rough isometries,
revealing more of the interest in this concept. The first of these is a
structure present in increasing rough isometries which we find quite
interesting, although, unfortunately, we have not found a way to use it
to our benefit in the sequel. The second of these is a slight variant
on rooted increasing rough isometries which will be much easier for us
to construct than general rough isometries; this variant is fundamental
to our construction in Section \ref{main_construction_section}.

\subsection{Increasing rough isometries as finite distributive
lattices} \label{increasing_ri_as_lattice_subsection}
In this subsection, we shall show that, given constants $(M,D,R)$ and
two finite subsets $A,B\subseteq\mathbb{N}\cup\{0\}$, both
containing 0, the
set of rooted increasing rough isometries from $A$ to $B$ with
constants $(M,D,R)$ is either empty or a finite distributive lattice.
This immediately implies a host of correlation inequalities (such as
the FKG inequality), as discussed below. However, although we consider
this to be a very interesting fact and a possibly useful structure, we
should mention at the outset that we do not use this fact in our
results and only include it here in the hope that it will prove useful
in further work on the problem.

We start with (see, e.g., \cite{AS00}, Chapter 6)
the following definition.
\begin{definition}
A finite partially ordered set $L$ is called a \textit{finite
distributive lattice} if any two elements $x,y\in L$ have a unique
minimal upper bound $x\vee y$ (called the \textit{join} of $x$ and $y$)
and a unique maximal lower bound $x\wedge y$ (called the \textit{meet}
of $x$ and $y$), such that, for any $x,y,z\in L$,
%
%
\begin{equation} \label{lattice_distributive_prop_eq}
x\wedge(y\vee z) = (x\wedge y)\vee(x\wedge z) .
\end{equation}
\end{definition}

Now, fix constants $(M,D,R)$ and finite subsets
$A,B\subseteq\mathbb{N}\cup\{0\}$, both containing~0, which are rooted
increasing r.i. with constants $(M,D,R)$. Let $L$ be the set of all
such rooted increasing r.i. mappings from $A$ to $B$. For $T_1,T_2\in
L$, we write $T_1\preceq T_2$ if, for all $x\in A$, we have $T_1(x)\le
T_2(x)$. We also define $T_1\vee T_2$ as $(T_1\vee T_2)\dvtx A\to B$,
$(T_1\vee T_2)(x):=\max(T_1(x),T_2(x))$ and, similarly, $(T_1\wedge
T_2)(x):=\min(T_1(x),T_2(x))$. It is clear that if $(T_1\vee T_2)\in
L$, then it is the unique minimal upper bound of $T_1$ and $T_2$ in $L$
and, similarly, that if $(T_1\wedge T_2)\in L$, then it is their unique
maximal lower bound. It is also clear that the distributive property
(\ref{lattice_distributive_prop_eq}) holds. Therefore, to show that $L$
is a finite distributive lattice, it remains to show the following.
\begin{lemma} \label{increasing_ri_as_lattice_lemma} For any
$T_1,T_2\in L$, we have $(T_1\vee T_2), (T_1\wedge T_2)\in L$. Or, in
words, the maximum and minimum of two rooted increasing r.i.'s with
constants $(M,D,R)$ are also rooted increasing r.i.'s with constants
$(M,D,R)$.
\end{lemma}

We remark that this lemma is not true for general rooted rough
isometries as it is easy to see, by means of examples, that the
monotonicity property is required.
\begin{pf*}{Proof of Lemma \protect\ref{increasing_ri_as_lattice_lemma}}
We shall show this for $T_1\vee T_2$, the proof for $(T_1\wedge T_2)$
being analogous (or even deducible from the $T_1\vee T_2$ case by
considering the reversed mappings). Letting $T:=T_1\vee T_2$, it is
clear that $T(0)=0$ and that $T$ is still (weakly) monotonic. We
continue by verifying property (ii) in the definition of r.i. (see
Figure \ref{distributive_lattice_prop_ii}). If we fix $b\in B$, then
there exist $x,y\in A$ with $|T_1(x) - b|\le R$ and $|T_2(y) - b|\le R$
and we may assume, without loss of generality, that $x\le y$. Of
course, if $T(x)=T_1(x)$, then property (ii) holds, hence we assume
that $T(x)=T_2(x)$. We obtain that $T_1(x)\le T(x)=T_2(x)\le T_2(y)$,
from which $|T(x) - b|\le R$ readily follows.

%
\begin{figure}

\includegraphics{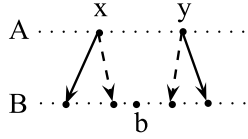}

\caption{$T_1$: solid line; $T_2$: dashed line.}
\label{distributive_lattice_prop_ii}
\end{figure}

Fixing $x,y\in A$, $x<y$, it remains to verify property (i) in the
definition of r.i. for $T$ and $x,y$ (see also Figure
\ref{distributive_lattice_prop_ii}). If $T(x)=T_i(x)$ and $T(y)=T_i(y)$
for $i=1$ or $i=2$, then the properties clearly hold since they hold
for $T_i$, hence we assume, without loss of generality, that
$T(x)=T_2(x)>T_1(x)$ and $T(y)=T_1(y)>T_2(y)$, to obtain
\[
\frac{1}{M}(y-x)-D\le T_2(y)-T_2(x)<T(y)-T(x)<T_1(y)-T_1(x)\le
M(y-x)+D ,
\]
proving the lemma.
\end{pf*}

The usefulness of the finite distributive lattice structure in
probability lies in the fact that it allows one to obtain correlation
inequalities in many cases. Following \cite{AS00}, Chapter 6, we have the following definition and theorem.
\begin{definition}
A probability measure $\mu$ on $L$ is called \textit{log-supermodular}
if, for all $T_1,T_2\in L$,
\[
\mu(T_1)\mu(T_2)\le\mu(T_1\vee T_2)\mu(T_1\wedge T_2) .
\]
\end{definition}

%
\begin{theorem}[(FKG inequality)] If $\mu$ is log-supermodular and
$f,g\dvtx L\to\mathbb{R}_+$ are increasing [in the sense that
$f(T_1)\le
f(T_2)$ whenever $T_1\preceq T_2$], then
\[
\mathbb{E}_\mu fg \le(\mathbb{E}_\mu f)(\mathbb{E}_\mu g) .
\]
\end{theorem}

In our case, one may take, for example, $\mu$ to be the uniform measure
on $L$, and, supposing $x,y\in A$, we may take $f(T_1)=T_1(x)$ and
$g(T_1)=T_1(y)$. We immediately obtain that when sampling a rough
isometry uniformly from $L$,the images of $x$ and $y$ are positively
correlated. This example may not be so impressive since the result is
intuitive, but, it is still not obvious how to prove this result
directly (for arbitrary r.i. $A$ and $B$) and the significant point is
that we obtained it here for free from the structure of $L$.

\subsection{Markov rough isometries}\label{Markov_ri_subsection}
In this subsection, we introduce a slightly different (but equivalent
up to constants) definition of a rooted increasing rough isometry which
will be much easier to work with in the sequel.
\begin{definition}
Two subsets $A,B\subseteq\mathbb{N}\cup\{0\}$ both containing $0$ are
\textit
{Markov rough isometric} if there exists a mapping $T\dvtx A\to B$ and
constants $M,F,R\ge0$ such that:
\begin{longlist}
\item $T(0)=0$;
\item if $x,y\in A$ and $x\ge y$, then $T(x)\ge T(y)$;
\item for all \textit{adjacent} $x,y\in A$ (i.e., with no point of $A$
between $x$ and $y$) with $T(x)\neq T(y)$, we have $\frac
{1}{M}|x-y|\le
|T(x)-T(y)|\le M|x-y|$;
\item for all $b\in T(A)$, we have $\max T^{-1}(b) - \min T^{-1}(b) \le F$;
\item for any $b\in B$, there exists $x\in A$ such that $|T(x)-b|\le R$.
\end{longlist}
\end{definition}

The reason for the name ``Markov rough isometry'' is that all of the
restrictions in the definition are, in some sense, local. To check that
a given mapping $T$ is a valid Markov rough isometry, one scans its
values on $A$ starting from $0$ and proceeding in increasing order. To
check the properties, one needs to remember the value of $T$ on a point
$x\in A$ only until one reaches a point $y>x$ with $T(y)>T(x)$ and, by
property (iv), this must happen after checking at most $F$ points.
Hence, there is a form of finite-memory property for Markov rough
isometries, which accounts for the name. Still, although they may
appear weaker at first, Markov rough isometries are equivalent to
rooted increasing rough isometries as follows.
\begin{lemma} \label{Markov_increasing_equiv_lemma}
Fix two subsets $A,B\subseteq\mathbb{N}\cup\{0\}$, both containing $0$.
\begin{enumerate}
\item If $T\dvtx A\to B$ is a Markov rough isometry with constants
$(M,F,R)$, then $T$ is a rooted increasing rough isometry with
constants $(2F+M,\frac{1}{2},R)$.
\item If $T\dvtx A\to B$ is a rooted increasing rough isometry with
constants $(M,D,R)$, then $T$ is a Markov rough isometry with constants
$(MD+M+D,MD,R)$.
\end{enumerate}
\end{lemma}
\begin{pf}

\begin{enumerate}
\item Let $T\dvtx A\to B$ be a Markov rough isometry with constants
$(M,F,R)$ and define $(\widetilde{M},\widetilde{D},\widetilde{R}):=(2F+M,\frac
{1}{2},R)$. To show that $T$ is a rooted increasing r.i. with constants
$(\widetilde{M},\widetilde{D},\widetilde{R})$, only property (i) in the definition
of rough isometry needs to be checked. If we let $x,y\in A$, $x<y$, and
first suppose that $T(x)\neq T(y)$, then we can find some $k\ge2$ and
a sequence of points of $A$, $x\le z^1_r<z^2_l\le z^2_r<\cdots
<z^{k-1}_l\le z^{k-1}_r<z^k_l\le y$, such that for each $i$, $z^i_r$ is
adjacent in $A$ to $z^{i+1}_l$, $T(z^i_r)\neq T(z^{i+1}_l)$,
$T(z^i_l)=T(z^i_r)$, $T(x)=T(z^1_r)$ and $T(y)=T(z^k_l)$ (Figure \ref
{Markov_increasing_equiv_pic} shows an example with $k=5$). Then
\begin{eqnarray*}
y-x &=& (y - z^k_l) + (z^k_l - z^{k-1}_r) + (z^{k-1}_r - z^{k-1}_l) +
\cdots+ (z^1_r - x) \\
&\le& kF + M\bigl(T(z^k_l)-z^{k-1}_r\bigr) + \cdots+ M\bigl(T(z^2_l) - T(z^1_r)\bigr) \\
&=& kF + M\bigl(T(y)-T(x)\bigr)
\end{eqnarray*}
and noting that $T(y)-T(x)\ge k-1$ [and, in particular, that
$T(y)-T(x)\ge1$], we obtain
\begin{eqnarray*}
y-x &\le& kF+M\bigl(T(y)-T(x)\bigr) \le2\bigl(T(y)-T(x)\bigr)F + M\bigl(T(y)-T(x)\bigr) \\
&=& \widetilde{M}\bigl(T(y)-T(x)\bigr) .
\end{eqnarray*}
The lower bound follows more easily:
\begin{eqnarray*}
y-x &\ge&(z^k_l - z^{k-1}_r) + (z^{k-1}_l - z^{k-2}_r) + \cdots+
(z^2_l - z^1_r) \\
&\ge&\frac{1}{M}\bigl(T(z^k_l - z^{k-1}_r) + \cdots+ T(z^2_l) - T(z^1_r)\bigr) =
\frac{1}{M}\bigl(T(y)-T(x)\bigr) .
\end{eqnarray*}

%
%
\begin{figure}

\includegraphics{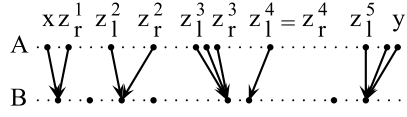}

\caption{An example for Lemma \protect\ref{Markov_increasing_equiv_lemma}
with $k=5$.}\label{Markov_increasing_equiv_pic}
\end{figure}

If we now suppose that $x,y\in A$, $x<y$, satisfy $T(x)=T(y)$, then
$y-\break x\le F$, hence
\[
0 = T(y)-T(x) \ge\frac{1}{2F+M}(y-x) - \frac{1}{2} = \frac
{1}{\widetilde
{M}}(y-x) - \widetilde{D}
\]
as required.
\item Let $T\dvtx A\to B$ be a rooted increasing rough isometry with
constants $(M,D,R)$ and define $(\widetilde{M},\widetilde{F},\widetilde
{R})=(MD+M+D,MD,R)$. To show that $T$ is a Markov r.i. with constants
$(\widetilde{M},\widetilde{D},\widetilde{R})$, only properties (iii) and (iv) in
the definition of Markov rough isometry need to be checked. If we let
$x,y\in A$ with $x$ adjacent to $y$ and $T(x)\neq T(y)$, then
\[
y-x\le M\bigl(T(y)-T(x)+D\bigr) \le(M+MD)\bigl(T(y)-T(x)\bigr) \le\widetilde{M}\bigl(T(y)-T(x)\bigr)
\]
and
\[
T(y)-T(x) \le M(y-x)+D \le(M+D)(y-x) \le\widetilde{M}(y-x) .
\]

If we now suppose that $x,y\in A$ satisfy $T(x)=T(y)$, then we have
\[
0=T(y)-T(x)\ge\frac{1}{M}(y-x)-D ,
\]
hence $y-x\le MD = \widetilde{F}$, as required.
\end{enumerate}
\upqed\end{pf}

We conclude this subsection by remarking that some properties of rooted
increasing rough isometries also hold for Markov rough isometries
(without the need to change the constants). First, it is trivial to
check the following (analogous to Proposition \ref
{increasing_RI_restriction_prop}).
\begin{proposition} \label{Markov_RI_restriction_prop}
If $A,B\subseteq\mathbb{R}$ are Markov rough isometric by a mapping
$T\dvtx
A\to B$ with constants $(M,F,R)$, then, for any $x,y\in A$ with $x<y$,
we have that $T$ restricted to $A\cap[x,y]$ is a Markov rough isometry
from $A\cap[x,y]$ to $B\cap[T(x),T(y)]$ with constants $(M,F,R)$.
\end{proposition}

Second,
we have the following proposition.
\begin{proposition}
Given $A,B\subseteq\mathbb{N}\cup\{0\}$, both containing 0, which
are Markov
rough isometric with constants $(M,F,R)$, the set $L$ of all Markov
rough isometries between them with constants $(M,F,R)$ is a finite
distributive lattice (with the same operations as defined in Section
\ref{increasing_ri_as_lattice_subsection}).
\end{proposition}

The proof is very similar to the proof of Lemma \ref
{increasing_ri_as_lattice_lemma} and is therefore omitted.

\section{Proof of equivalence theorems} \label
{equivalence_thms_proofs_section}
We start with the proof of Proposition \ref{zero_one_event_prop}.
\begin{pf*}{Proof of Proposition \protect\ref{zero_one_event_prop}}
We will use the well-known fact that a Poisson process on $\mathbb{R}$
with the
shift operation on $\mathbb{R}$ is ergodic. We also note that the
event $E$
that $A$ and $B$ are rough isometric with constants $(M,D,R)$ is
measurable with respect to $A$ and $B$. Next, we note that for any
fixed realization of $B$, the event $E_B$ that $A$ is rough isometric
to $B$ with constants $(M,D,R)$ is translation invariant (with respect
to translations of $A$), hence, by ergodicity, it has probability $0$
or $1$. Analogously, for any fixed realization of $A$, the event $E_A$
that $A$ is rough isometric to $B$ with constants $(M,D,R)$ is also
translation invariant (with respect to translations of $B$) and hence
has probability $0$ or $1$. It now follows from the independence of $A$
and $B$ that $E$ itself has probability $0$ or $1$.
\end{pf*}

We continue with the proof of Proposition \ref{Poisson_perc_equiv_prop}.
\begin{pf*}{Proof of Proposition \protect\ref{Poisson_perc_equiv_prop}}
$\mbox{(ii)}\to\mbox{(i)}$. This is trivial.

$\mbox{(i)}\to\mbox{(ii)}$. Suppose that claim (i) holds for some $\alpha, \beta>0$.
Fix $\gamma>0$ and consider two Poisson processes $A$ and $C$, with
intensities $\alpha$ and $\gamma$, respectively. Note that they can be
coupled by first sampling $A$ and then letting the points of $C$ be
$\{\frac{\alpha}{\gamma}x | x\in A\}$. Now, observe that under
this coupling, $A$ and $C$ are r.i. with constants $(\frac{\alpha
}{\gamma}, 0, 0)$ under the trivial mapping $T\dvtx A\to C$ defined by
$T(x):=\frac{\alpha}{\gamma} x$.

In the same way, if we fix some $\delta>0$, then we can couple two
Poisson processes $B$ and $D$, with intensities $\beta$ and $\delta$,
respectively, so that they are rough isometric a.s. Considering now two
such independent Poisson processes $A$ and $B$, and the processes $C$
and $D$ which are coupled to them, we find that $C$ and $D$ are also
independent and that they are rough isometric a.s. by transitivity of
the rough isometry relation since $C$ and $A$ are rough isometric a.s.
by our coupling, $A$ and $B$ are rough isometric a.s. using (i) and $B$
and $D$ are rough isometric a.s. by our coupling.

By means of similar transitivity arguments, to prove that (iii) and
(iv) are equivalent to (i) and (ii), it is enough to establish that for
any $\alpha>0$ and $0<p<1$, a Poisson process $A$ of intensity $\alpha$
and a Bernoulli percolation $\mathcaligr{A}$ with parameter $p$ can be
coupled to be rough isometric a.s. We now show this. If we fix $\alpha$
and $p$ to have a coupling first sample $A$, then $\mathcaligr{A}$ will
have a point at the integer $n$ if and only if $A$ has at least one
point in the interval $[nc, (n+1)c)$, where $c=-\frac{\log
(1-p)}{\alpha
}$ is chosen so that this is indeed a coupling. Now, define a mapping
$T\dvtx\mathcaligr{A}\to A$ by $T(n):=x_n$, where $x_n$ is some point of
$A$ in the interval $[nc, (n+1)c)$, say the smallest one. It is easy to
see that $T$ is a rough isometry with constants $(\max(c,\frac{1}{c}
),c,c)$ since if $n,n+k\in\mathcaligr{A}$, then $(k-1)c\le
T(n+k)-T(n)\le(k+1)c$.
\end{pf*}

\subsection{Proof of Theorem \protect\ref{discrete_to_oriented_equiv_thm}}
We first prove Lemma \ref{cut_point_lemma}.
\begin{pf*}{Proof of Lemma \protect\ref{cut_point_lemma}}
Let $A$ and $B$ be two independent Bernoulli percolations on $\mathbb{Z}$.
First, by  assumption, there exist constants $(M,D,R)$ such that $A$
and $B$ are rough isometric a.s. with constants $(M,D,R)$. Let $\Omega
^1$ be this event.

Second, for $k,l,m\in\mathbb{Z}$ with $k<l<m$, let $\Omega
^2_{k,l,m}$ be the
event that $k,l,m\in A$, $l$ and $m$ are adjacent in $A$ (no point of
$A$ lies between them) and $m-l\ge\frac{l-k}{M^2}-\frac{2D}{M}$.
Noting that for fixed $k$, we have $\mathbb{P}(\Omega^2_{k,l,m})\le
2^{-(m-l)}1_{(m-l)\ge c(l-k)-C}$ for some $C,c>0$, we get
\[
\sum_{(m,l | m>l>k)} \mathbb{P}(\Omega^2_{k,l,m}) \le\sum_{(l | l>k)}
2^{-c(l-k)+C+1} < \infty.
\]
The Borel--Cantelli lemma then implies that with probability 1, only
finitely many $\Omega^2_{k,l,m}$ occur for a fixed $k$.

Third, we condition on the events $\Omega^1$ and the event that for
each $k$, only finitely many $\Omega^2_{k,l,m}$ occur. We fix two
realizations $A$ and $B$, and let $T\dvtx A\to B$ be the r.i. between
them. We will show (a deterministic claim) that there exists a cut
point for $T$. To see this, fix $a\in A$ and let $b:=T(a)\in B$, noting
that if there are only finitely many $u_n\in A$ with $u_n>a$ and
$T(u_n)<b$, then if we take $x$ to be the largest of these $u_n$, $x$
will satisfy ($\alpha$) in the definition of cut point. Analogously, if
there were only finitely many $v_n\in A$ with $v_n>a$ and $T(v_n)>b$,
then ($\beta$) (in the definition of cut point) would be satisfied for
some $x$. Hence, we assume, by way of contradiction, that there are
infinitely many such $u_n$ and $v_n$. Since only finitely many $x\in A$
can be mapped to $b$, we must have infinitely many pairs $v,u\in A$,
adjacent in $A$ with $a<v<u$, $T(v)>b$ and $T(u)<b$. Each such pair
must satisfy
\[
u-v\ge\frac{1}{M}\bigl(T(v)-T(u) - D\bigr) \ge\frac{1}{M}\bigl(T(v) - b - D\bigr) \ge
\frac{1}{M}\biggl(\frac{1}{M}(v-a)-2D\biggr) ,
\]
but this is a contradiction since only finitely many $\Omega
^2_{a,l,m}$ occur.
\end{pf*}
\begin{pf*}{Proof of Theorem \protect\ref{discrete_to_oriented_equiv_thm}}
$\mbox{(ii)}\to\mbox{(i)}$. Let $A$ and $B$ be two independent Bernoulli percolations
on $\mathbb{Z}$. With probability $\frac{1}{4}$, they both contain $0$.
Conditional on this event, let $(A^+,0)$ be the rooted Bernoulli
percolation on $\mathbb{N}$ obtained from $A$ by considering only the
nonnegative integers. Define similarly the independent $(A^-,0)$
obtained by considering the nonpositive integers, and the independent
$(B^+,0)$ and $(B^-,0)$. By (ii), there exist constants $(M,D,R)$ such
that, with positive probability, $(A^+,0)$ is rooted r.i. to $(B^+,0)$
and $(A^-,0)$ is rooted r.i. to $(B^-,0)$ with these constants. Denote
these r.i. mappings by $T^+$ and $T^-$, respectively. Let the map
$T\dvtx A\to B$ be the map whose restriction to $A^+$ is $T^+$ and whose
restriction to $A^-$ is $T^-$. It is then easy to check directly from
the definition that $T$ is a r.i. of $A$ to $B$ with constants
$(M,2D,R)$. This shows that, with positive probability, $A$ and $B$ are
r.i., but according to Propositions \ref{zero_one_event_prop} and
\ref{Poisson_perc_equiv_prop}, $A$ and $B$ are r.i. with probability
$0$ or $1$. Hence, $A$ and $B$ are rough isometric a.s.

$\mbox{(i)}\to\mbox{(ii)}$. Let $p$ be the probability that two independent rooted
Bernoulli percolations on $\mathbb{N}$ are rooted r.i. We need to show that
$p>0$. Let $A$ and $B$ be two independent Bernoulli percolations on
$\mathbb{Z}$. For $n,m\in\mathbb{Z}$, let $A_n^+$ be all points of
$A$ not smaller
than $n$ and let $A_n^-$ be all points of $A$ not larger than $n$;
similarly define $B_m^+$ and $B_m^-$. Let $E^{+}_{n,m}$ be the event
that $n\in A$, $m\in B$ and there exists a rooted r.i. between
$(A_n^+,n)$ and $(B_m^+,m)$; similarly define $E^{-}_{n,m}$ using
$A_n^+$ and~$B_m^-$. Note that $\mathbb{P}(E^{+}_{n,m}) = \mathbb
{P}(E^{-}_{n,m}) =
\frac{p}{4}$. Now, since by (i) and Lemma \ref{cut_point_lemma}, with
probability 1, there exists a r.i. $T\dvtx A\to B$ with a cut point
$x\in
A$, we get that $\mathbb{P}(\bigcup_{n,m} (E^{+}_{n,m}\cup
E^{-}_{n,m})) = 1$.
This implies that $p>0$, proving the claim.
\end{pf*}

\subsection{Proof of Theorem \protect\ref
{finite_and_infinite_perc_equiv_thm} and related lemmas}
We start with the following proof.
\begin{pf*}{Proof of Lemma \protect\ref{big_gap_deterministic_lemma}}
Let $z\in A$ be the largest point such that $T(z)\le T(x)$. Note that
$z$ must be finite [since $T(0)=0$ and $A$ is infinite] and that $z\ge
y>x$. First, note that for large enough $L$ (as a function of $M$ and
$D$),
%
%
\begin{equation} \label{use_of_L_ineq}
z-x\ge y-x\ge\frac{T(x)-T(y)-D}{M} \ge\frac{L-D}{M} \ge\frac
{L}{2M} .
\end{equation}
Second, let $w:=\operatorname{Succ}(z)$. Note that, by definition of
$z$, we have
$T(w)>T(x)\ge T(z)$, hence,
\begin{eqnarray*}
w-z&\ge&\frac{T(w)-T(z)-D}{M} > \frac{T(x)-T(z)-D}{M}\\
& \ge&\frac
{1}{M}\biggl(\frac{z-x}{M}-2D\biggr)\\
&=& \frac{z-x}{M^2}-\frac{2D}{M}
\end{eqnarray*}
and, by combining this inequality with (\ref{use_of_L_ineq}), we see
that if $L$ is large enough (as a function of $M$ and $D$), then
$w-z\ge
\frac{z-x}{2M^2}$, as required.
\end{pf*}

We next show the following.
\begin{pf*}{Proof of Lemma \protect\ref{big_gap_prob_estimate_lemma}}
For any fixed $z\in\mathbb{N}$, $\mathbb{P}(\operatorname
{Gap}(z)\ge k)\le2^{-(k-1)}$ (with equality
if $k$ is a positive integer). Hence, by a union bound,
\begin{eqnarray*}
\mathbb{P}(E^w_{L,M})
&\le&
\biggl(\frac{L}{2M}+1\biggr)2^{-({L}/({4M^3})-1)} +
\sum
_{i=\lceil{L}/({2M})\rceil}^\infty2^{-({i-1})/({2M^2})}\\
&\le&
C\biggl(\frac
{L}{M}+1\biggr)e^{-c{L}/({M^3})}.
\end{eqnarray*}
\upqed\end{pf*}

We continue with the following proof.
\begin{pf*}{Proof of Lemma \protect\ref{rooted_ri_restriction_lemma}}
Letting $x_i\in A$ be the $i$th point of $A$ and $a_i$ be the $i$th
point of $B$, we choose $m$ so that $a_m=\max_{1\le i\le n} T(x_i)$
[i.e., the minimal $m$ such that $T(A(n))\subseteq B(m)]$. First, for
any $x,y\in A(n)$, we have
\[
\frac{1}{M}|x-y| - D\le|T(y)-T(x)| \le M|x-y| + D ,
\]
by the properties of $T$. Second, to reach a contradiction, assume
that for some $b\in B(m)$ and for all $x\in A(n)$, $|T(x)-b|>L$. Since
$T$ is a rooted r.i. with constants $(M,D,R)$, there must exist some
$y\in A$, $y>x_n$ with $|T(y)-b|\le R$; furthermore, by the minimality
of $m$, there must exist some $x\in A(n)$ with $T(x)>b+L$, hence $x\le
x_n<y$ and $T(x)-T(y)\ge L-R$. By Lemma
\ref{big_gap_deterministic_lemma}, there exists some $z\in A$, $z\ge y$
and $z-x\ge\frac{L-R}{2M}$ such that $\operatorname{Gap}(z)\ge\frac
{z-x}{2M^2}$.
But, then, in particular, $z>x_n$ and $\operatorname{Gap}(z)\ge
\max(\frac{L-R}{4M^3}, \frac{z-x_n}{2M^2} )$, which
contradicts the fact that $E^{x_n}_{L-R,M}$ does not hold for $A$.
\end{pf*}

Finally, we have the following proof.

\begin{pf*}{Proof of Theorem \protect\ref{finite_and_infinite_perc_equiv_thm}}
$\mbox{(i)}\to\mbox{(ii)}$. Let $(A,0)$ and $(B,0)$ be two independent rooted
Bernoulli percolations on $\mathbb{N}$ and let $E$ be the event that
they are
rooted r.i. with constants $(M,D,R)$. Suppose that $\mathbb{P}(E)\ge
r$ for
some $r>0$. On the event~$E$, let $T\dvtx A\to B$ be such a rooted r.i. Fix
$n\ge1$, let $x_n\in A$ be the $n$th point of~$A$, fix $L>R$ and let
$E^{x_n}_{L-R,M}$ be the event from Lemma
\ref{big_gap_prob_estimate_lemma}. Note that since $x_n$ is a stopping
time for the percolation $A$ (i.e., $\{x_n>k\}$ only depends on whether
$i\in A$ for $0\le i\le k$) and since $E^{x}_{L,M}$ only depends on the
future of $x$ (i.e., on the events $\{i\in A\}_{i>x}$), we have, by
Lemma \ref{big_gap_prob_estimate_lemma}, that $\mathbb
{P}(E^{x_n}_{L-R,M})\le
C(\frac{L-R}{M}+1)e^{-c({L-R})/({M^3})}$ for some absolute constants
$C,c>0$. Hence, for each fixed $0<p<r$, we can choose $L$ sufficiently
large (uniformly in~$n$) so that $\mathbb{P}(E\cap
(E^{x_n}_{L-R,M})^c)\ge p$;
we fix such a pair of $p$ and $L$. We are now done since, on the event
$E\cap(E^{x_n}_{L-R,M})^c$, Lemma~\ref{rooted_ri_restriction_lemma}
gives that $T$ restricted to $A(n)$ is a rooted r.i. of $A(n)$ to some
initial segment of $B$ with constants $(M,D,L)$.

$\mbox{(ii)}\to\mbox{(i)}$. Let $E_n$ be the event that $A(n)$ is rooted r.i. to some
initial segment of $B$ with constants $(M,D,R)$, so that by  assumption that
$\mathbb{P}(E_n)\ge p>0$ for all~$n$. Let $E:=\limsup E_n$, so that,
by Fatou's
lemma, $\mathbb{P}(E)\ge\limsup\mathbb{P}(E_n)\ge p$. Let $(A,B)\in
E$, that is, $A$
and $B$ are two realizations of rooted Bernoulli percolation on
$\mathbb{N}$
such that, for an infinite sequence $n_k\to\infty$ (depending on $A$
and $B$), there exists a rooted r.i. $T_{n_k}$ from $A(n_k)$ to some
initial segment of $B$ with constants $(M,D,R)$. We now deduce that $A$
and $B$ are themselves rooted r.i. with constants $(M,D,R)$. Let $x_i$
be the $i$th point of $A$. To define $T\dvtx A\to B$, we need to pick
$a_i\in B$ such that $T(x_i):=a_i$; we do this by induction. Since
$x_1=0$, we also choose $a_1:=0$. Assume that we have already chosen
$\{a_i\}_{i=1}^{N-1}$ for some $N\ge2$ in such a way that there exists
an infinite sequence $n^j:=n_{k_j}$ such that $T_{n^j}$ agrees with $T$
on $\{x_i\}_{i=1}^{N-1}$. To choose $a_N$, we note that
$\{T_{n^j}(x_N)\}_j$ is a finite set since, for example, for each $j$,
$\frac{x_N}{M}-D\le T_{n^j}(x_N)\le Mx_N + D$. Hence, we can choose
$a_N$ in such a way that it agrees with an infinite subsequence of
$\{T_{n^j}\}_j$. In this way, we obtain $T$.

To see that $T$ is a rooted r.i. with
constants $(M,D,R)$, we note that for each $x,y\in A$, by our
construction, there exists some $k$ such that $T_{n_k}$ agrees with $T$
on $x$ and $y$. Hence, $\frac{|x-y|}{M} - D\le|T(x) - T(y)|\le M|x-y|
+ D$. Next, we fix $b\in B$, choose $N$ sufficiently large that
$x_N\ge M(b+R+D)$ and choose $k$ so that $T_{n_k}$ agrees with $T$ on
$\{x_i\}_{i=1}^N$. Since $T_{n_k}$ is a rooted r.i., there exists $x\in
A$ such that $|T(x)-b|\le R$. We cannot have $x>X_N$ since, otherwise,
$|T(x)|\ge\frac{x}{M}-D>\frac{x_N}{M}-D\ge b+R$. Hence, $x\le x_N$ and
so $|T(x)-b|\le R$, as required. This completes the proof of the
theorem.
\end{pf*}

\section{The main construction} \label{main_construction_section}
In this section, we shall prove Theorem \ref{main_thm}. Let us recall
the setting. We are given two independent rooted Bernoulli percolations
$(A,0)$ and $(B,0)$ on $\mathbb{N}$. We will show that, for any large enough
$n$ (independent of $A$ and $B$), there exists a Markov rough isometry
from $A(n)$ to some initial segment of $B$ with constants $(10\sqrt
{\log
_2 n},10\sqrt{\log_2 n},10\sqrt{\log_2 n})$ and with probability
$1-2^{-8\sqrt{\log_2 n}}$. As explained before, existence of a Markov
rough isometry is a stronger statement than existence of a general
rough isometry since Markov rough isometries are monotone and, by Lemma
\ref{Markov_increasing_equiv_lemma}, the same mapping\vspace*{1pt} will also be a
rooted increasing rough isometry with constants $(30\sqrt{\log_2 n},
\frac{1}{2}, 10\sqrt{\log_2 n})$. The reason we construct a Markov
rough isometry rather than an increasing rooted rough isometry is that
we will frequently rely on the fact that one can check the validity of
a Markov rough isometry by simply looking at local configurations (as
explained in Section \ref{Markov_ri_subsection}).

We fix $n$ very large. It would be convenient for us to assume that
$M,F$ and $R$ are integers, hence we choose $0.99<\alpha<1$ (depending
on $n$) so that $\alpha\sqrt{\log_2 n}$ is an integer. We then let
$M=F=R:=10\alpha\sqrt{\log_2 n}$. We also\vspace*{-1pt} introduce a new parameter,
$K:=2^{\alpha\sqrt{\log_2 n}}= (2^M )^{{1}/{10}}$, whose use will
be made clear in the sequel.

Given a sorted sequence $U:=(0,x_1,x_2, \ldots, x_L)\subseteq\mathbb
{N}\cup
\{0\}
$ (where we allow $L$ to be infinite), we define some notation. For a
point $t\in U$, let $s^U(t)$ or, equivalently, $s_1^U(t)$ be its
successor point in $U$; similarly, let $s_k^U(t)$ be its $k$th
successor point in $U$ and define $s_0^U(t):=t$. We call the
quantity
$g^U(t):=s^U(t)-t$ the \textit{gap at $t$}. When the set $U$ is clear
from the context, we sometimes omit the superscript and simply write
$s_k(t)$ and $g(t)$.

We will sometimes refer to $U$ equivalently by its \textit{gap
sequence} $\{G^U(i)\}_{i=1}^{L}$, defined by $G^U(i):=x_i-x_{i-1}$.

Let $A$ and $B$ be two independent rooted Bernoulli percolations
$(A,0)$ and $(B,0)$ on $\mathbb{N}$. Note that for $A$ and $B$, the sequences
$G^A$ and $G^B$ are simply i.i.d. $\operatorname{Geom}(\frac{1}{2})$
random variables.

We shall call a gap \textit{short} if it less than or equal to $M$,
otherwise we call it \textit{long}.

\subsection{Partitioning into blocks}
\begin{figure}[b]

\includegraphics{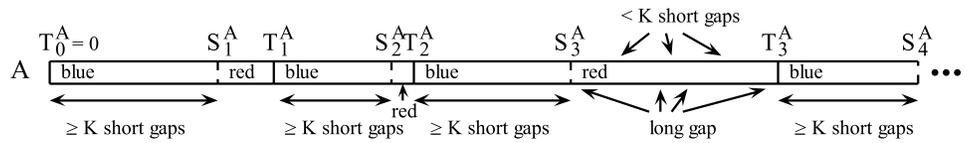}

\caption{A sample of the first three blocks followed by the blue
segment of the fourth block. The third red segment has long and short
gaps indicated.}\label{division_into_blocks_with_T_and_S}
\end{figure}

The first thing we will do is to partition $A$ and $B$ into blocks
(which overlap at their end points). Let us first describe this
partition informally and then give a rigorous definition. Each block
will consist of two parts, a ``blue'' initial segment followed by a
``red'' segment. A blue segment is a segment of the percolation points
containing only short gaps (of length $\le M$). A red segment is a
segment of the percolation points starting with a long gap (of length
$>M$) and ending just before $K$ short gaps (see Figure \ref
{division_into_blocks_with_T_and_S}).

%

More formally, to define blocks in $A$, we define a sequence of times
inductively, $T_0^A:=0$, and, for each $k\ge1$,
%
%
\begin{eqnarray}
S_k^A &:=& \min\{t\in A | t>T^A_{k-1}, g(t)>M\} ,\nonumber\\[-8pt]\\[-8pt]
T_k^A &:=& \min\{t\in A | t>S_k^A, g(s_i(t))\le M \mbox{ for all $0\le
i\le K-1$}\} .\nonumber
\end{eqnarray}

For each $k\ge1$, $S_k^A$ is the first point in $A$ after $T_{k-1}^A$
and immediately preceding a gap longer than $M$, and $T_k^A$ is the
first point in $A$ after $S_k^A$ which precedes $K$ short gaps.

The points of $A$ in the segment $[T_{k-1}^A,T_k^A]$ constitute the
$k$th \textit{block of $A$}. In each block, the \textit{blue segment}
consists of the points in $[T_{k-1}^A,S_k^A]$. By definition (except
possibly for the first block), the blue segment has at least $K$ short
gaps (and no long gaps). It is followed by a \textit{red segment},
consisting of the points in $[S_k^A,T_k^A]$, which starts with a long
gap and continues until the starting point of a run of $K$ short gaps
(not including that run). Note that the red segment may contain many
long gaps or as few as one. Also, it must start with a long gap and end
immediately after a long gap. The first block is different from the
rest since it may have less than $K$ gaps in its blue segment. However, letting
%
%
\begin{equation} \label{k_short_gaps_event_def}
E_0^A:=\{\mbox{$A$ starts with at least $K$ short gaps}\} ,
\end{equation}
we have $\mathbb{P}(E_0^A) = (1-\frac{1}{2^M})^K \ge1-\frac{K}{2^M} =
1-2^{-9\alpha\sqrt{\log_2 n}}$. We emphasize that, conditioned on
$E_0^A$, the distribution of blocks after subtracting their starting
points (or, equivalently, when looking at their gap sequences) is
i.i.d. and we shall refer to that common distribution as $\mathcaligr
{L}^{\mathrm
{block}}$, or, in words, the distribution of a \textit{rooted block}.

We partition $B$ in the same way, into blocks analogously defining
$T_k^B$, $S_k^B$ and~$E_0^B$.

It will be useful to define the distributions of blocks and of blue and
red segments precisely, which we now proceed to do.
\begin{definition}
We say that $X\sim\operatorname{Geom}_{\le M}(\frac{1}{2})$ if $X$
is distributed
like a $\operatorname{Geom}(\frac{1}{2})$ random variable conditioned
to be less than
or equal to $M$. We say that $Y\sim\operatorname{Geom}_{>M}(\frac
{1}{2})$ if $Y$ is
distributed like a $\operatorname{Geom}(\frac{1}{2})$ random variable
conditioned to
be larger than $M$ or, in other words, as $M+\operatorname{Geom}(\frac
{1}{2})$.
\end{definition}

The following observation will be useful in the sequel. It is also true
in much greater generality.
\begin{lemma} \label{stochastic_dominance_observ_lemma}
The $\operatorname{Geom}_{\le M}(\frac{1}{2})$ distribution is stochastically
dominated by the $\operatorname{Geom}(\frac{1}{2})$ distribution.
\end{lemma}
\begin{pf}
Define a coupling of $(X,Y)$ with $X\sim\operatorname{Geom}_{\le
M}(\frac{1}{2})$
and $Y\sim\operatorname{Geom}(\frac{1}{2})$ using the following
algorithm: take an
infinite sequence $(Z_i)_{i=1}^{\infty}$ of i.i.d. $\operatorname
{Geom}(\frac{1}{2})$
random variables, and let $Y=Z_1$ and $X=Z_i$, with $i$ the minimal
index for which $Z_i\le M$. It is then clear that $X\le Y$ a.s.
\end{pf}
\begin{definition}
For a given integer $L>0$, say that $U:=(0,x_1,x_2,\ldots,\break
x_L)\subseteq
\mathbb{N}\cup\{0\}$ is distributed $\mathcaligr{L}_L^{\mathrm
{blue}}$, or in words,
distributed as a \textit{rooted blue segment of length $L$} if
$(x_i-x_{i-1})_{i=1}^L$ are i.i.d. $\operatorname{Geom}_{\le M}(\frac
{1}{2})$ (where $x_0:=0$).
\end{definition}
\begin{lemma}\label{block_structure_lemma}
Let $B=(0,x_1,\ldots, x_P, x_{P+1},\ldots, x_Q)$ be a rooted block,
with $U:=(0,x_1,\ldots, x_P)$ being its blue segment and
$(x_P,x_{P+1},\ldots,x_Q)$ being its red segment. Also, let
$V:=(0,x_{P+1}-x_P,\ldots, x_Q-x_P)=(0,y_1,\ldots, y_{Q-P})$ be the red
segment minus its starting point. Then:
\begin{enumerate}
\item$U$ and $V$ are independent;
\item$U$ is distributed $\mathcaligr{L}_P^{\mathrm{blue}}$, where $P$
is a random
variable distributed $\operatorname{Geom}(\frac{1}{2^M})-1$,
conditioned to be at
least $K$ [or, in other words, $P\sim K-1+\operatorname{Geom}(\frac
{1}{2^M})$],
independently of the lengths of the gaps in the block;
\item the distribution of $V$ is characterized by:
\begin{enumerate}[(a)]
\item[(a)] $y_1\sim\operatorname{Geom}_{>M}(\frac{1}{2})$, independently
of the other gaps;
\item[(b)] $y_2,\ldots,y_{Q-P}$ are the concatenation of $N\sim
\operatorname{Geom}
((1-\frac{1}{2^M})^K)-1$ subsequences which are i.i.d., given $N$. Each
such subsequence starts with $Z$ gaps, each having distribution
$\operatorname{Geom}
_{\le M}(\frac{1}{2})$ independently of each other and where $Z$ is
distributed $\operatorname{Geom}(\frac{1}{2^M})-1$, conditioned to be
less than $K$.
The subsequence then continues with one last gap distributed
$\operatorname{Geom}
_{>M}(\frac{1}{2})$ independently of the other gaps.
\end{enumerate}
\end{enumerate}
\end{lemma}
\begin{pf}
The red segment begins at the first long gap of a block. It is
clear that knowing the lengths of all of the gaps previous to this gap
does not provide any additional information on the length of this or
the following gaps.

The first, say, blue segment of $A$ contains all of the gaps up
to the first long gap from the beginning of $A$. The length of this run
of short gaps is $\operatorname{Geom}(\frac{1}{2^M})-1$ and it is
independent of the
lengths of the short gaps in it. Hence,  conditioned that this run
of short gaps contains at least $K$ gaps, we obtain the
characterization given in the lemma.

The first, say, red segment of $A$ is defined to start where the
first run of short gaps of $A$ ends and to continue until just before a
run of at least $K$ short gaps. Hence, it can be described in the
following way. First, since it ends a run of short gaps, it has to
start with a long gap. Since the gaps in $A$ are i.i.d. and all we know
about this gap is that it is long, its size will be independent of the
size of all other gaps [but distributed $\operatorname
{Geom}_{>M}(\frac{1}{2})]$. We
then test to see if the following $K$ gaps are all short. If they are,
then we end the red segment, otherwise we include the run of short gaps
coming afterward and the long gap following it in the red segment. We
now continue in the same manner with another independent trial to see
if the next $K$ gaps are all short. If so, we end, otherwise we include
them and the long gap at their end in the red segment. These
independent trials continue until we finally find a run of at least $K$
short gaps. Hence, the number of trials is geometric (but we subtract
one since once we succeed, we do not concatenate anything to the red
segment) and its success parameter is $(1-\frac{1}{2^M})^K$, which
is\vspace*{1pt}
the probability of seeing $K$ short gaps in a row. When a trial fails,
it means that the number of short gaps after it is less than $K$.
Since, a priori, the number of short gaps is $\operatorname
{Geom}(\frac{1}{2^M})-1$,
we have that $Z$, the number of short gaps following a failed trial, is
$\operatorname{Geom}(\frac{1}{2^M})-1$, conditioned to be less than
$K$. Finally, the
lengths of the short gaps themselves are unaffected by the number of
short gaps in a run, hence they are all $\operatorname{Geom}_{\le
M}(\frac{1}{2})$,
independently of everything else. Similarly, the length of the long gap
which ends a run of short gaps is $\operatorname{Geom}_{>M}(\frac{1}{2})$,
independently of everything else.
\end{pf}

\begin{definition}
We say that a vector having the distribution of the vector $V$ of the
previous lemma is distributed $\mathcaligr{L}^{\mathrm{red}}$ or, in words,
distributed as a \textit{rooted red segment}.
\end{definition}

\subsection{Properties of blocks}
In this subsection, we will prove some basic properties of rooted blue
and red segments which will be useful for our construction in the
sequel. We start with two properties of red segments.
\begin{lemma} \label{red_segment_prop_lemma}
Let $V\sim\mathcaligr{L}^{\mathrm{red}}$, $X$ be the number of long
gaps in $V$ and
$\{b_i\}_{i=1}^X$ be their lengths. There then exist $\beta,\gamma>0$
such that
\begin{eqnarray*}
\mathbb{P}\biggl(X>\frac{1}{8}\sqrt{\log_2 n}\biggr)&=&o \biggl(\frac{1}{n^{1+\beta
}} \biggr),\\
\mathbb{P}\Biggl(\sum_{i=1}^X b_i \ge3\log_2 n\Biggr)&=& o \biggl(\frac
{1}{n^{1+\gamma}} \biggr).
\end{eqnarray*}
\end{lemma}
\begin{pf}
By Lemma \ref{block_structure_lemma}, we know that $X\sim
\operatorname{Geom}
((1-\frac
{1}{2^M})^K)$. Hence,
%
%
\begin{eqnarray}\label{number_of_red_gaps_estimate}
\mathbb{P}\biggl(X> \frac{1}{8}\sqrt{\log_2 n}\biggr) &=& \biggl[1-\biggl(1-\frac{1}{2^M}\biggr)^K
\biggr]^{
{1}/{8}\sqrt{\log_2 n}} \le\biggl(\frac{K}{2^M} \biggr)^{{1}/{8}\sqrt{\log_2
n}} \nonumber\\[-8pt]\\[-8pt]
&=& 2^{-{9}/{8}\alpha\log_2 n} =
o \biggl(\frac{1}{n^{1+\beta}} \biggr)\nonumber
\end{eqnarray}
for some $\beta>0$, proving the first claim. Now, conditioned on $X$,
the $\{b_i\}_{i=1}^X$ are i.i.d. with distribution $\operatorname
{Geom}_{>M}(\frac
{1}{2})$, that is, with distribution $M + \operatorname{Geom}(\frac
{1}{2})$. Hence,
\begin{eqnarray*}
\mathbb{P}\Biggl(\sum_{i=1}^X b_i\ge
3\log_2 n | X\Biggr)
&=&
\sum_{s\ge3\log_2 n}\mathop{\mathop{\sum}_{b_1+\cdots+b_X=s}}_{b_i>M} 2^{-\sum_{i=1}^X
(b_i-M)}\\
&=&2^{MX}\sum_{s\ge3\log_2 n} 2^{-s}\#\{b_1+\cdots+ b_X=s | b_i>M\}
\\
&\le& 2^{10\alpha\sqrt{\log_2 n} X}\sum_{s\ge3\log_2 n} 2^{-s}s^{X},
\end{eqnarray*}
so, denoting $E:=\{X\le\frac{1}{8}\sqrt{\log_2 n}\}$, we have, for
large enough $n$ and some $C>0$,
\begin{eqnarray*}
\mathbb{P}\Biggl(\Biggl\{\sum_{i=1}^X b_i\ge3\log_2 n\Biggr\}\cap E\Biggr)
&\le&
2^{{10}/{8}\alpha
\log_2 n}\sum_{s\ge3\log_2 n} 2^{-s}s^{{1}/{8}\sqrt{\log_2
n}}
\\
&\le& 2^{{5}/{4}\alpha\log_2 n}\sum_{s\ge3\log_2 n}
2^{-s}s^{{1}/{8}\sqrt{{s}/{3}}} \\
&\le& 2^{{5}/{4}\alpha\log_2 n}\sum_{s\ge3\log_2 n}
2^{-{4}/{5}s} \\
&\le& C2^{{5}/{4}\alpha\log_2 n - {12}/{5}\log_2 n} = o
\biggl(\frac{1}{n^{1+\widetilde{\gamma}}} \biggr)
\end{eqnarray*}
for some $\widetilde{\gamma}>0$. Hence, by (\ref
{number_of_red_gaps_estimate}), we have $\mathbb{P}(\sum_{i=1}^X
b_i\ge3\log_2
n) \le o (\frac{1}{n^{1+\beta}} ) + o (\frac{1}{n^{1+\widetilde{\gamma}}}
)$, proving the second claim.
\end{pf}

We continue with three properties of blue segments. We start with the
following, simple, lemma.
\begin{lemma} \label{stopping_time_prop_for_blue_lemma}
For a given integer $L>0$ and $U:=(0,x_1,x_2,\ldots, x_L)\sim\mathcaligr{L}
_L^{\mathrm{blue}}$, if $0\le T\le L$ is a stopping time in the sense
that the event $\{T\le k\}$ depends only on $\{x_i\}_{i=1}^k$, then,
conditioned on $T$, on the event $\{T<L\}$, the partial rooted segment
$V:=(0,x_{T+1}-x_T,\ldots, x_L-x_T)$ is distributed $\mathcaligr{L}
_{L-T}^{\mathrm{blue}}$.
\end{lemma}
\begin{pf}
Consider the gap sequence $G^U=(x_1, x_2-x_1, \ldots, x_L-x_{L-1})$. By
definition, its elements are i.i.d. $\operatorname{Geom}_{\le M}(\frac
{1}{2})$. If we
let $A_k:=\{T=k\}$ for $k<L$ and $B$ be an event that depends only on
$(x_{k+1}-x_k, \ldots, x_L-x_k)$, then, since $A_k$ is determined by
$(x_1, \ldots, x_k)$ and these are, in turn, determined by $(x_1,
x_2-x_1, x_k-x_{k-1})$, we have that $A_k$ and $B$ are independent.
Hence, conditioned on $A_k$, the probability of $B$ remains the same,
implying that $(x_{k+1} - x_k, \ldots, x_L-x_k)$ are still i.i.d.
$\operatorname{Geom}
_{\le M}(\frac{1}{2})$, proving the claim.
\end{pf}
\begin{lemma} \label{division_to_subsegments_lemma}
Fix integers $L,Z>0$ and let $U:=(0,x_1,\ldots, x_L)\sim\mathcaligr{L}
_L^{\mathrm
{blue}}$. Divide the points of $U$ into subsegments according to the
following algorithm: the first subsegment consists of $(0,x_1, x_2,
\ldots, x_{l_1})$ with $l_1$ maximal such that $x_{l_1}\le Z$; by
induction for $i\ge2$, the $i$th subsegment consists of
$(x_{l_{i-1}+1},\ldots, x_{l_i})$ with $l_i$ maximal such that
$x_{l_i}-x_{l_{i-1}+1}\le Z$. Let $Y$ be the number of subsegments
required to cover all $L$ points. We claim that
\[
\mathbb{P}\biggl(Y>\frac{3L}{Z}\biggr)\le e^{-cL}
\]
for some $c>0$.
\end{lemma}
\begin{pf}
First, note that the event $Y>m$ is contained in the event $x_L>mZ$.
Hence, $\mathbb{P}(Y>\frac{3L}{Z})\le\mathbb{P}(\sum_{i=1}^L
G_i>3L)$, where the $G_i$
are i.i.d. $\operatorname{Geom}_{\le M}(\frac{1}{2})$ random
variables. Since a
$\operatorname{Geom}
_{\le M}(\frac{1}{2})$ random variable is stochastically dominated by a
$\operatorname{Geom}(\frac{1}{2})$ random variable, by Lemma \ref
{stochastic_dominance_observ_lemma}, we get, by standard large
deviation estimates, that $\mathbb{P}(\sum_{i=1}^L G_i>3L)\le
e^{-cL}$ for some
$c>0$, as claimed.
\end{pf}

The following lemma is a major ingredient in our rough isometry construction.
\begin{lemma} \label{waiting_time_for_red_match}
Let $(G_i)_{i=1}^\infty$ be i.i.d. $\operatorname{Geom}_{\le M}(\frac
{1}{2})$ random
variables, termed \textit{gaps}. Let $m>0$, $M>a_1,\ldots, a_m>0$ and
$d_1,\ldots, d_{m-1}\ge0$ be given integers. We consider the $\{a_i\}$
as representing minimal required gap lengths and the $\{d_i\}$ as
representing inter-gap distances. Say that a position $l$ is \textit
{valid} if $G_l\ge a_1$, $G_{l+d_1+1}\ge a_2$, $G_{l+d_1+1+d_2+1}\ge
a_3, \ldots, G_{l+m-1+\sum_{i=1}^{m-1} d_i}\ge a_m$. If we let $Z$ be
the minimal valid position, then, for any $a>0$ and $s:=\sum_{i=1}^m a_i$,
\[
\mathbb{P}(Z> \lceil a2^s\rceil)\le e^{-{a}/{m^2}}.
\]
\end{lemma}
\begin{pf}
Let $I_l$ be the event that $l$ is a valid position. Then
\[
\mathbb{P}(I_l) = \prod_{i=1}^m \biggl(\frac{1}{2^{a_i-1}}-\frac{1}{2^M}
\biggr)\bigg/\biggl(1-\frac
{1}{2^M}\biggr) \ge\prod_{i=1}^m \frac{1}{2^{a_i}} = \frac{1}{2^s} .
\]
For a given position $l$, let us denote by $C_l:=(l,l+d_1+1,\ldots,
l+m-1+\sum_{i=1}^{m-1} d_i)$ the \textit{comb at position $l$} and say
that two positions $l,k$ overlap if their combs intersect, that is, if
$C_l\cap C_k\neq\varnothing$ (see Figure \ref{combs_on_gap_sequence}).
Note that if $F\subseteq\mathbb{N}$ is a subset of positions, no two
of which
overlap, then $\{I_l\}_{l\in F}$ are independent.

%
\begin{figure}

\includegraphics{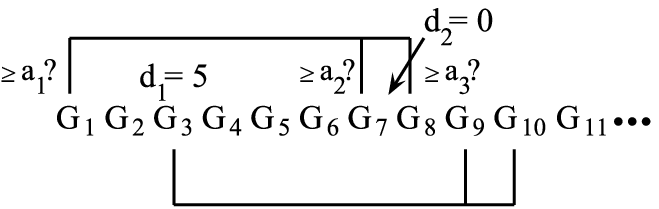}

\caption{A comb at two nonoverlapping positions.}\label{combs_on_gap_sequence}
\end{figure}

Fixing an integer $N>0$, to bound $\mathbb{P}(Z>N)$, we wish to choose
a large
collection of positions $F\subseteq\{1,\ldots,N\}$, no two of which
overlap. We note that a given comb $C_l$ may only intersect at most
$m(m-1)$ other combs $C_k$ since each overlapping position $k$ uniquely
determines a pair of coordinates $1\le i,j\le m$, $i\neq j$, such that
the $i$th coordinate of $C_l$ is equal to the $j$th coordinate of $C_k$
by, say, the smallest element of $C_l\cap C_k$. Hence, we can find such
a collection $F$ with, say, $|F|\ge\lceil\frac{N}{m^2}\rceil$, by
means of a greedy algorithm. Thus, we obtain the bound
\begin{eqnarray*}
\mathbb{P}(Z>N)
&\le&
\mathbb{P}\biggl(\bigcap_{l\in F} I_l^c\biggr) = \prod_{l\in F} \bigl(1-\mathbb{P}
(I_l)\bigr) \\
&\le& \biggl(1-\frac{1}{2^s}\biggr)^{\lceil{N}/{m^2}\rceil} \le
e^{-2^{-s}\lceil
{N}/{m^2}\rceil}
\end{eqnarray*}
and the claim follows by taking $N:=\lceil a2^s\rceil$.
\end{pf}
\begin{remark} \label{stopping_time_valid_pos_remark}
We point out that in the notation of the previous lemma, the position
$Z+m-1+\sum_{i=1}^{m-1} d_i$ is a stopping time for the process $\{
G_i\}
_{i=1}^\infty$.
\end{remark}

\subsection{The construction}
A major part of the construction of the rough isometry between $A$ and
$B$ will be constructing a rough isometry between a block of $A$ and
the beginning of a blue segment of $B$ or, alternatively, constructing
a rough isometry between the beginning of a blue segment of $A$ and a
block of $B$. The following theorem gives conditions under which this
is possible with high probability.
\begin{theorem} \label{block_mapping_theorem}
Fix integers\vspace*{1pt} $L_1, L_2$ satisfying $L_2\ge\max(K,L_1)$ and \mbox{$L_1\ge
\frac{K}{2}$}. Let $U^1:=(0,x^1_1,\ldots, x^1_{L_1})\sim\mathcaligr
{L}_{L_1}^{\mathrm
{blue}}$, $U^2:=(0,x^2_1,\ldots, x^2_{L_2})\sim\mathcaligr
{L}_{L_2}^{\mathrm
{blue}}$ and $V\dvtx(0,y_1,\break\ldots, y_N)\sim\mathcaligr{L}^{\mathrm
{red}}$, where
$N$ is random, with $U^1, U^2, V$ independent. Construct the segment
$W:=(0,x^1_1,\ldots,x^1_{L_1}, x^1_{L_1}+y_1, \ldots, x^1_{L_1}+y_N)$
by concatenating $U^1$ and~$V$. There then exists a random integer
$1\le S\le L_2$ which is a stopping time for $U^2$ conditioned on $W$.
That is, the event $\{S\le l\}$ is measurable with respect to $W$ and
$\{x_i^2\}_{i=1}^l$, satisfying the conditions that if $E=\{S\le\max
(\frac{K}{2},\frac{L_1}{\sqrt{\log_2 n}} )\}$, then:
\begin{longlist}
\item$\mathbb{P}(E)=1-o (\frac{1}{n^{1+\delta}} )$ for some $\delta>0$;
\item on the event $E$, there exists a Markov rough isometry $T_1$ from
$W$ to $U^2\cap[0,x_S^2]$ with constants $(M,F,R)$ such that the last
point of $W$ is mapped to $x_S^2$ and it is the only point mapped to $x_S^2$;
\item on the event $E$, there exists a Markov rough isometry $T_2$ from
$U^2\cap[0,x_S^2]$ to $W$ with constants $(M,F,R)$ such that $x_S^2$
is mapped to the last point of $W$ and it is the only point mapped to
the last point of $W$.
\end{longlist}
\end{theorem}

Let us show how to prove Theorem \ref{main_thm} using Theorem \ref
{block_mapping_theorem}. We first require the following definition.
\begin{definition}
For a given integer $L\ge0$, a sorted infinite sequence
$U:=(x_1,x_2,x_3,\ldots)\subseteq\mathbb{N}\cup\{0\}$ is said to be
distributed
as a \textit{Bernoulli percolation with $L$ initial short gaps} if the
rooted sequence $V:=(0,x_2-x_1,x_3-x_1,\ldots)$ is distributed as a
rooted Bernoulli percolation on $\mathbb{N}$ conditioned to have its
first $L$
gaps short and its next gap long.
\end{definition}

The proof is by induction: for each stage $0\le j\le n$, we shall have
an event $E_j$ denoting whether or not the $j$th stage was successful,
with $\mathbb{P}(E_j | \{E_i\}_{i=0}^{j-1})=1-o (\frac{1}{n^{1+\delta
}} )$ for
$j\ge1$ and $\mathbb{P}(E_0)\ge1-2^{-9\alpha\sqrt{\log_2 n}+1}$.
Conditioned
on $\bigcap_{i=0}^j E_i$, the following random variables are defined:
%
\begin{longlist}
\item two positions $P_j^A\in A$ and $P_j^B\in B$, with $P_j^A\ge s_j^A(0)$;
\item a Markov rough isometry $T_j\dvtx A\cap[0,P_j^A]\to B\cap
[0,P_j^B]$ with constants $(M,F,R)$ satisfying $T_j(P_j^A)=P_j^B$, with
$P_j^A$ being the only source of $P_j^B$;
\item two numbers $L^A_j$ and $L^B_j$, with $\max(L_j^A,L_j^B)\ge K$
and $\min(L_j^A,L_j^B)\ge\frac{K}{2}$.
\end{longlist}
Also, conditioned on all of these random variables, the distribution of
$A\cap[P_j^A,\infty)$ is that of a Bernoulli percolation with $L^A_j$
initial short gaps and, independently, the distribution of $B\cap
[P_j^B,\infty)$ is that of a Bernoulli percolation with $L^B_j$ initial
short gaps.

This implies Theorem \ref{main_thm} since if all events $\{
E_j\}_{j=0}^n$ occur, then $T_n\dvtx A\cap[0,P_n^A]\to B\cap[0,P_n^B]$
is a Markov rough isometry with constants $(M,F,R)$ and $P_n^A\ge
s_n^A(0)$, and hence, by Proposition \ref{Markov_RI_restriction_prop},
we know that its restriction to the first $n$ points of $A$ is a Markov
r.i. to some initial segment of $B$ with constants $(M,F,R)$, as the
theorem requires. The probability that $\{E_j\}_{j=0}^n$ occur is at
least $(1-2^{-9\alpha\sqrt{\log_2 n}+1})(1-o (\frac{1}{n^{1+\delta}})
)^n > 1-2^{-8\sqrt{\log_2 n}}$ for large enough $n$, as required.

%
\begin{figure}

\includegraphics{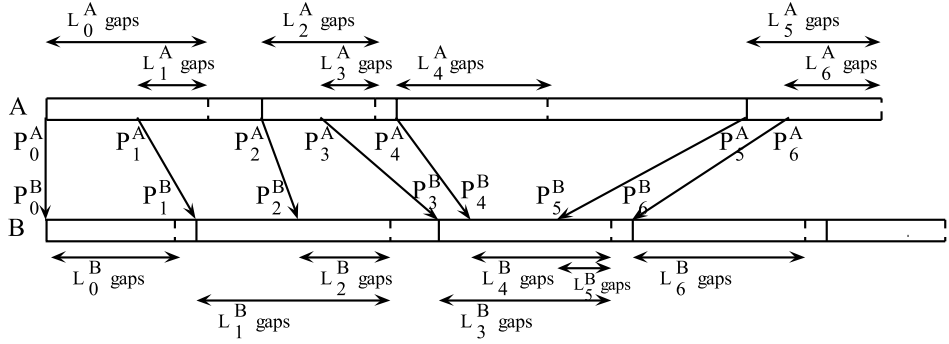}

\caption{Illustration of the mapping of the first percolation into the
other using the induction procedure. The blue and red segments of
blocks are depicted as in Figure
\protect\ref{division_into_blocks_with_T_and_S}. In this example,
$L_0^A\ge
L_0^B, L_1^A\le L_1^B, \ldots.$}\label{block_squeezed_into_blue_pic}
\end{figure}

Let us show the above induction (see Figure \ref
{block_squeezed_into_blue_pic}). For $j=0$, the event $E_0:=E_0^A\cap
E_0^B$ [recall (\ref{k_short_gaps_event_def})], $P_j^A=P_j^B=0$, $T_0$
is just defined on $0\in A$ by $T_0(0)=0$ and, on the event $E_0$, we
let $L_0^A$ be the length of the first blue segment of
$A$ and $L_0^B$ be the length of the first blue segment of $B$.
It is easy
to see that all of the properties stated above hold.

Now, suppose that $\{E_i\}_{i=0}^{j-1}$ have occurred and that we have
already constructed the above random variables up to stage $j-1$ with
the above properties. We condition on $\bigcap_{i=0}^{j-1} E_i,
P_{j-1}^A, P_{j-1}^B, T_{j-1}, L_{j-1}^A$ and $L_{j-1}^B$. There are
two cases to consider:
\begin{enumerate}
\item$L^B_{j-1}\ge L^A_{j-1}$ [note that this also implies
$L^B_{j-1}\ge K$, by property (iii) above]. Let
$Q_j^A:=s_{L^A_{j-1}}(P_{j-1}^A)$, $Q_j^B:=s_{L^B_{j-1}}(P_{j-1}^B)$.
By the induction assumption, the segment $A\cap[P_{j-1}^A,Q_j^A]$
translated to start at $0$ is distributed $\mathcaligr
{L}_{L^A_{j-1}}^{\mathrm
{blue}}$ and the segment $B\cap[P_{j-1}^B,Q_j^B]$ translated to start
at $0$ is distributed $\mathcaligr{L}_{L^B_{j-1}}^{\mathrm{blue}}$. Let $P_j^A$
denote the end of the red segment which follows $A\cap
[P_{j-1}^A,Q_j^A]$, that is,
\[
P_j^A:=\min\{x\in A | x>Q_j^A, g(s_i(x))\le M\mbox{ for all $0\le
i\le
K-1$}\} ,
\]
and let $L^A_j$ be the number of short gaps of $A$ after $P_j^A$, that is,
\[
L^A_j:=\max\{N | g(s_i(P_j^A))\le M\mbox{ for all $0\le i\le N-1$}\} .
\]
Note that, by definition of $P_j^A$, we have $L^A_j\ge K$. We now
invoke Theorem \ref{block_mapping_theorem} with the following
parameters: $U^1$ is the segment $A\cap[P_{j-1}^A,Q_j^A]$ translated to
start at $0$, $V$ is the segment $A\cap[Q_j^A, P_j^A]$ translated to
start at $0$ and $U^2$ is the segment $B\cap[P_{j-1}^B,Q_j^B]$
translated to\vspace*{2pt} start at $0$ ($W$ is then $A\cap[P_{j-1}^A,P_j^A]$
translated to start at $0$). The theorem gives us $S$, which is a
stopping time for $U^2$ conditioned on $U^1$ and $V$. Let $E_j$ be the
event $E$ of that theorem, that is,
\[
E_j:=\biggl\{S\le\max\biggl(\frac{K}{2},\frac{L^A_{j-1}}{\sqrt{\log_2 n}} \biggr)\biggr\}.
\]
According to part (ii) of that theorem, on the event $E_j$, we have a
Markov rough isometry $\widetilde{T}_j\dvtx W\to U^2\cap[0,s_S(0)]$ with
constants $(M,F,R)$. Let $P_j^B:=s_S(P_{j-1}^B)$ and $L_j^B:=L_{j-1}^B
- S$, and note that, on the event $E_j$, $L_j^B\ge\max
(L^A_{j-1},K)-S\ge\frac{1}{2}\max(L^A_{j-1},K) \ge\frac{K}{2}$.
Finally, to construct $T_j$ we ``concatenate'' $T_{j-1}$ and $\widetilde
{T}_j$, that is,
\[
T_j(x):=\cases{ T_{j-1}(x),&\quad $x\in A, x\le P_{j-1}^A$,\cr
\widetilde{T}_j(x-P_{j-1}^A)+P_{j-1}^B,&\quad $x\in A, P_{j-1}^A\le x\le
P_j^A$.}
\]
Note that $T_j$ is indeed a Markov rough isometry with constants
$(M,F,R)$ since $T_{j-1}$ and $\widetilde{T}_j$ are, and since there is a
unique preimage to $P_{j-1}^B$. Also note that by Lemma \ref
{stopping_time_prop_for_blue_lemma}, we have that, conditioned on $E_j$
and $S$, the distribution of $B\cap[P_j^B,\infty)$ is that of a
Bernoulli percolation with $L^B_j$ initial short gaps. Hence, $E_j$,
$P_j^A$, $P_j^B$, $T_j$, $L_j^A$ and $L_j^B$ satisfy the requirements
of the induction step.
\item$L^A_{j-1}\ge L^B_{j-1}$. The induction step in this case is
performed in the same way as in the first case, but with the roles of
$A$ and $B$ interchanged and using part (iii) of Theorem \ref
{block_mapping_theorem} instead of part (ii).
\end{enumerate}

All that remains is to prove Theorem \ref{block_mapping_theorem}, which
we now do.
\begin{pf*}{Proof of Theorem \protect\ref{block_mapping_theorem}}
We divide the proof into several parts:
\begin{enumerate}
\item First, consider $U^1$. Applying the algorithm of Lemma
\ref{division_to_subsegments_lemma} to $U^1$ with $Z=F$, we obtain a
division of $U^1$ into $Y$ subsegments. Denote these by $U^1_1,\ldots,
U^1_Y$. If we let $\Omega_1:=\{Y\le\frac{3L_1}{F}\}$, then, by Lemma
\ref{division_to_subsegments_lemma}, there exists a $c>0$ such that
\[
\mathbb{P}(\Omega_1^c) \le e^{-cL_1} \le e^{-c{K}/{2}}.
\]

\item Now, consider $V$. Let $X$ be the number of long gaps in $V$ and
let $\{b_i\}_{i=1}^X$ be their lengths. Let
$\Omega_2:=\{X\le\frac{1}{8}\sqrt{\log_2 n}\}$ and
$\Omega_3:=\{\sum_{i=1}^X b_i\le3\log_2 n\}$. Then, by Lemma
\ref{red_segment_prop_lemma}, for some $\beta,\gamma>0$,
%
%
\begin{eqnarray}\label{red_segment_estimates}
\mathbb{P}(\Omega_2^c)&\le& o \biggl(\frac{1}{n^{1+\beta}}
\biggr),\nonumber\\[-8pt]\\[-8pt]
\mathbb{P}(\Omega_3^c)&\le& o \biggl(\frac{1}{n^{1+\gamma}} \biggr).\nonumber
\end{eqnarray}

\item We continue to consider $V$. Let $(z_i)_{i=1}^X\subseteq V$ be
the starting points of the long gaps in $V$ ($z_1=0$), that is,
$g(z_i)>M$ for all $i$. They divide $V$ into $X-1$ subsegments,
$V^1:=V\cap[s(z_1),z_2], \ldots, V^{X-1}:=V\cap[s(z_{X-1}),z_X]$. By
the (structure) Lemma \ref{block_structure_lemma}, we know that each
subsegment conditioned on its length and translated to start at $0$ is
distributed as a blue segment of that length. Let us again employ the
algorithm of Lemma \ref{division_to_subsegments_lemma} with $Z=F$ to
each of these subsegments to divide them further into
``sub-subsegments.'' Let $Y_i$ be the number of sub-subsegments in the
division of $V^i$ and denote them by $(V^i_j)_{j=1}^{Y_i}$ (for each
$j$, $V^i_j\subseteq V^i$). Let $\Omega_4:=\{\sum_{i=1}^{X-1}
Y_i\le\frac{K}{20}\}$. To bound $\sum_{i=1}^{X-1} Y_i$, we consider the
blue segment $\widetilde{V}$ obtained from $V$ by deleting all of its
long gaps and translating to start at $0$. More precisely, write
$V^i=(y^i_0,\ldots, y^i_{N^i})$, where $N_i$ is the number of gaps in
$V^i$, let $\widetilde{V}^i:=(0,y^i_1-y^i_0,\ldots,
y^i_{N^i}-y^i_0)=:(0,\widetilde{y}^i_1,\ldots, \widetilde{y}^i_{N^i})$
and then define
\begin{eqnarray*}
\widetilde{V}&:=&\Biggl(\underbrace{0,\widetilde{y}^1_1,\ldots,
\widetilde
{y}^1_{N^1}}_{\widetilde{V}^1},
\underbrace{\widetilde{y}^2_1+\widetilde{y}^1_{N^1},\ldots,
\widetilde
{y}^2_{N^2}+\widetilde{y}^1_{N^1}}_{\mathrm{Translated}\ \widetilde
{V}^2},\ldots, \\
&&\hspace*{35pt}\underbrace{\widetilde{y}^{X-1}_1+\sum_{j=1}^{X-2} \widetilde
{y}^j_{N^j}, \ldots,
\widetilde{y}^{X-1}_{N^{X-1}} + \sum_{j=1}^{X-2} \widetilde
{y}^j_{N^j}}_{\mathrm{Translated}\ \widetilde{V}^{X-1}}\Biggr) .
\end{eqnarray*}
$\widetilde{V}$ is a blue segment as a concatenation of many
independent blue segments. We also apply the algorithm of Lemma \ref
{division_to_subsegments_lemma} to $\widetilde{V}$ with $Z=F$ to divide
it into $\widetilde{Y}$ subsegments. It is clear from the algorithm
that $\widetilde{Y}\le\sum_{i=1}^{X-1} Y_i$, but since, in the passage
from $V$ to $\widetilde{V}$, we only removed $X$ long gaps, one must
also check that
%
%
\begin{equation} \label{sum_of_Ys_bound}
\sum_{i=1}^{X-1} Y_i\le\widetilde{Y}+X .
\end{equation}
We recall that $N$ is the number of gaps in $V$ and note that by the
(structure) Lemma \ref{block_structure_lemma}, $N\le1+K(X-1)\le KX$.
We wish to show that
%
%
\begin{equation} \label{tilde_Y_prob_estimate}
\mathbb{P}\biggl(\widetilde{Y}>\frac{K}{21}\biggr)= o \biggl(\frac{1}{n^{1+\beta}} \biggr).
\end{equation}
For this, we divide the problem into three cases:
\begin{itemize}
\item$N>\frac{KF}{70}$. This implies that $X>\frac{F}{70}>\frac
{1}{8}\sqrt{\log_2 n}$, which we know, by (\ref
{red_segment_estimates}), to have probability at most $o (\frac
{1}{n^{1+\beta}} )$.
\item$\frac{KF}{70}\ge N>\frac{K}{21}$. Applying Lemma \ref
{division_to_subsegments_lemma}, we have
\begin{eqnarray*}
\mathbb{P}\biggl(\widetilde{Y}>\frac{K}{21}, \frac{KF}{70}\ge N>\frac
{K}{21}\biggr)&\le& \mathbb{P}
\biggl(\widetilde{Y}>\frac{3N}{F}, \frac{KF}{70}\ge N>\frac{K}{21}\biggr)\\
&\le&\mathbb{E} e^{-cN}\mathbf{1}_{ ({KF}/{70}\ge N>{K}/{21} )}\le
e^{-c{K}/{21}}.
\end{eqnarray*}
\item$N\le\frac{K}{21}$. On this event, we certainly must have
$\widetilde{Y}\le\frac{K}{21}$,
\end{itemize}
and (\ref{tilde_Y_prob_estimate}) follows. Using (\ref{red_segment_estimates}), (\ref
{sum_of_Ys_bound}) and (\ref{tilde_Y_prob_estimate}), we deduce, for
large enough $n$, that
\[
\mathbb{P}(\Omega_4^c)\le\mathbb{P}\biggl(\widetilde{Y}>\frac{K}{21}\biggr)+\mathbb{P}(\Omega_2^c) =
o \biggl(\frac
{1}{n^{1+\beta}} \biggr).
\]
\item We now consider the gap sequence $G:=G^{U^2}=(g_i)_{i=1}^{L_2}$
of $U^2$ and the sequence $\widetilde{G}:=(g_i)_{i=Y}^{\infty}$ ($Y$
was defined in the first item of the proof), where we have extended the
sequence to be infinite by concatenating an i.i.d. sequence
$(g_i)_{i=L_2+1}^{\infty}$ of $\operatorname{Geom}_{\le M}(\frac
{1}{2})$ random
variables, independent of everything else. We apply Lemma \ref
{waiting_time_for_red_match} to $\widetilde{G}$ with the parameters
$m=X, a_i:= \lceil\frac{b_i}{M}\rceil$ (recall that $\{b_i\}_{i=1}^X$
are the lengths of the long gaps of $V$) and $d_i=Y_i-1$, to obtain
$Z$, the first valid position along $\widetilde{G}$ (``valid position''
was defined in the lemma). Let $\Omega_5:=\{Z\le K^{{3}/{4}}\}$.
Let $s:=\sum_{i=1}^X a_i$ and choose $a:=2^{{1}/{4}\sqrt{\log_2
n}}$. Then, by the lemma,
\[
\mathbb{P}(Z> \lceil a2^s\rceil| X, \{b_i\}_{i=1}^X, \{Y_i\}
_{i=1}^{X-1})\le
e^{-{a}/{X^2}}.
\]
Since, on the events $\Omega_2$ and $\Omega_3$, we have $s\le\frac
{3\log_2 n}{M} + X\le\frac{9}{20}\sqrt{\log_2 n}$, we obtain, for
large enough $n$,
\begin{eqnarray*}
&&\mathbb{P}(Z> K^{{3}/{4}}, \Omega_2, \Omega_3 | X, \{b_i\}
_{i=1}^X, \{
Y_i\}_{i=1}^{X-1}) \\
&&\qquad \le\mathbb{P}(Z> \lceil a2^s\rceil, \Omega_2, \Omega_3 | X, \{
b_i\}
_{i=1}^X, \{Y_i\}_{i=1}^{X-1}) \le e^{-{a}/{X^2}}.
\end{eqnarray*}
Hence, $\mathbb{P}(\Omega_5^c) = o (\frac{1}{n^{1+\beta}} ) + o
(\frac
{1}{n^{1+\gamma}} )$.\vspace*{1pt}

\item Finally, we construct the required time $S$, event $E$ and Markov
rough isometries $T_1$ and $T_2$. We define
\[
S:=Y+Z+X+\sum_{i=1}^{X-1} (Y_i-1)=Y+Z+1+\sum_{i=1}^{X-1} Y_i .
\]
We note that, just as in Remark \ref{stopping_time_valid_pos_remark},
conditioned on $W$ [in particular, on $Y,X$ and $(Y_i)_{i=1}^{X-1}$],
the time $S$ is a stopping time for $U^2$. We define the event
\[
\widetilde{E}:=\Omega_1\cap\Omega_2\cap\Omega_4\cap\Omega_5 .
\]
Note that, by the previous calculations, $\mathbb{P}(\widetilde{E}^c)=o
(\frac
{1}{n^{1+\delta}} )$ for some $\delta>0$. On the event $\widetilde{E}$,
we have
\[
S\le\frac{3L_1}{F}+K^{{3}/{4}}+1+\frac{K}{20}\le\max\biggl(\frac
{K}{2},\frac{L_1}{\sqrt{\log_2 n}} \biggr),
\]
hence the event $E$ of the theorem satisfies $E\supseteq\widetilde{E}$.
On the event $E$, we now construct $T_1$ (see Figure \ref
{mapping_of_block_to_blue_segment_pic}). $T_2$ is constructed
%
%
\begin{figure}

\includegraphics{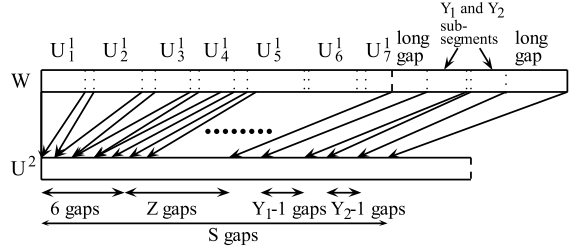}

\caption{Illustration of the constructed rough isometry. In the
picture, $Y=7$ and $X=3$. When mapping $U^1$ to $U^2$, we start mapping
points one-to-one rather than many-to-one, starting from subsegment
$j_0:=4$.}\label{mapping_of_block_to_blue_segment_pic}
\end{figure}
analogously, using the fact that $R=F$. First, we define $T_1$ on the
points of $U^1$ in such a way that $T_1(x^1_{L_1}) = x^2_{Y+Z}$. Note
that, on the event $E$,
%
%
\begin{eqnarray}\label{L_1_Y_and_Z_estimate}
L_1 - Y &=& L_1-S+Z+1+\sum_{i=1}^{X-1} Y_i \nonumber\\[-8pt]\\[-8pt]
&\ge& L_1 - \max\biggl(\frac
{K}{2},\frac{L_1}{\sqrt{\log_2 n}} \biggr) + Z\ge Z\nonumber
\end{eqnarray}
since $L_1\ge\frac{K}{2}$. We start mapping the points of $U^1$ to
$U^2$ according to the subsegment $U^1_j$ that the point of $U^1$ is
in. More precisely, we consider all of the points of $U^1$ in order and
if a point $x^1_i\in U^1_j$, then we define $T(x^1_i):=x^2_{j-1}$
(where $x^2_0$ is defined to be $0$). By the definition of the
subsegments $(U^1_j)_{j=1}^Y$, for all $j$, we will have $\max
T^{-1}(x^2_{j-1}) - \min T^{-1}(x^2_{j-1})\le F$, as required.
Furthermore, since $x^1_{i+1} - x^1_i\le M$ for all $i$ and
$x^2_{j+1}-x^2_j\ge1$, we will not expand or contract any distance by
more than $M$. We stop mapping in this way when we reach a point
$x^1_{i_0}$ belonging to $U^1_{j_0}$ which satisfies $L_1-i_0 =
Y+Z-(j_0-1)$. Such a point must be reached for some $0\le i_0\le L_1$,
by (\ref{L_1_Y_and_Z_estimate}). From this point on, we map the
points\vspace*{1pt}
sequentially as follows: $T(x^1_{i_0+l}):= x^2_{j_0-1+l}$ for $0\le
l\le L_1-i_0$. As before, no distances are expanded or contracted by
more than $M$.

We continue to define $T_1$ at the points of $W$ which follow
$x^1_{L_1}$ (the translated points of $V$). Recalling that
$(z_i)_{i=1}^X\subseteq V$ are the starting points of the long gaps in
$V$ (the corresponding points of $W$ are $z_i+x^1_{L_1}$), we
construct\vspace*{1pt}
the remainder of the mapping $T_1$ by induction on $1\le j\le X$. Note
that since $z_1=0$, we have already defined $T_1(z_1+x^1_{L_1}) =
x^2_{Y+Z}$. Define further $T_1(s(z_1+x^1_{L_1})):= x^2_{Y+Z+1}$; this
is the $j=1$ stage. Note that by the definition of $Z$, we did not
contract the gap of $W$ by more than $M$ (and, of course, we did not
expand it since we mapped to a short gap).

For $2\le j\le X$, let $R_j:=\sum_{i=1}^{j-1} Y_i$. Suppose that we
have already constructed the mapping $T_1$ to be a Markov rough
isometry with constants $(M,F,R)$ from $W\cap[0,s(z_{j-1}+x^1_{L_1})]$
to $U^2\cap[0,x^2_{Y+Z+1+R_{j-1}}]$ in such a way that
$T_1(s(z_{j-1}+x^1_{L_1}))=x^2_{Y+Z+1+R_{j-1}}$ and that it is the only
source of $x^2_{Y+Z+1+R_{j-1}}$. Recall that we have divided the
subsegment $V^{j-1}$ into sub-subsegments $(V^{j-1}_k)_{k=1}^{Y_{j-1}}$
and consider a point $x_l\in
W\cap[s(z_{j-1}+x^1_{L_1}),z_j+x^1_{L_1}]$. There then exists some
$k$ such that $x_l-x^1_{L_1}\in V^{j-1}_k$. Define
$T_1(x_l):=x^2_{Y+Z+1+R_{j-1} + (k-1)}$. Note that this is consistent
with the definition of $s(z_{j-1}+x^1_{L_1})$, that
$T_1(z_j+x^1_{L_1})=x^2_{Y+Z+R_j}$, that by the choice of the
sub-subsegments for each point $x^2_l\in
U^2\cap[x^2_{Y+Z+1+R_{j-1}},x^2_{Y+Z+R_j}]$, we have $\max
T_1^{-1}(x^2_l) - \min T_1^{-1}(x^2_l)\le F$, and that since we are
mapping gaps not larger than $M$ to gaps of size between $1$ and $M$,
no distance was expanded or contracted by more than $M$ (whenever two
points are mapped to different images). Finally, define
$T_1(s(z_j+x^1_{L_1})):= x^2_{Y+Z+1+R_j}$. Again, by the choice of $Z$,
this mapping did not contract the gap of $W$ by more than $M$ (and, of
course, we did not expand it since we mapped to a short gap).
Continuing this procedure until $j=X$ completes the construction of the
map $T_1\dvtx W\to U^2\cap[0,x^2_S]$, as required.
\end{enumerate}
\upqed\end{pf*}

\section*{Acknowledgments}
I would like to thank Itai Benjamini for
introducing me to the problem and encouraging me to work on it, and Ori
Gurel-Gurevich and Gady Kozma for suggesting ways to attack the problem
similar to the ones I have here employed. I would also like to thank
Gideon Amir, Steve Evans and Yuval Peres for useful conversations and
discussions concerning this problem, and Mikl\'os Ab\'ert for useful
discussions on the history of the problem and its general version.
Finally, I wish to thank Guillaume Obozinsky and Nicholas Crawford for
noticing and correcting an error in a previous version of Figure
\ref{monotone_non_monotone_L_example}.

%

%
\printaddresses

\end{document}